%% file: submission_file.tex
\documentclass[hidelinks,onefignum,onetabnum]{siamart250211}
\input{MMS_shared}

\usepackage{verbatim}
\usepackage{latexsym, graphicx, epsfig}
\usepackage{pgfplots}
\usepackage{mathrsfs}
\usepackage{color}
\usepackage{array}
\usepackage{amsfonts,amsmath,amssymb}
\usepackage{float}
\usepackage{multirow}
\usepackage{booktabs}
\usepackage{subcaption}
\usepackage{tikz}
\usepackage{diagbox}
\usepackage[toc,page]{appendix}
\usepackage{algorithm}
\usepackage[numbers,sort&compress]{natbib}
\usepackage{enumitem}  
\usepackage{xr}
\usepackage{mathtools}
\usepackage{mathabx}
\usepackage{comment}
\usepackage{cancel}
\usepackage{ulem}
\usepackage{cleveref}

\newcommand{\e}{\varepsilon}

\newsiamthm{assumption}{Assumption}
\numberwithin{equation}{section}
\numberwithin{figure}{section}

\title{A Bi-fidelity numerical method for velocity discretization of Boltzmann equations  
\thanks{
    Submitted to arXiv in July, 2025.
    \funding{
    N. Crouseilles thanks L. Liu for her kind hospitality during his visit to Chinese University of Hong Kong. He also acknowledges the EUROfusion Consortium, funded by the European Union via the Euratom Research and Training Programme, grant agreement No 101052200 EUROfusion(views and opinions expressed are however those of the authors only and do not necessarily reflect those of the European Union or the European Commission). 
    Z. Hao acknowledges the computational resources provided by The Chinese University of Hong Kong during his employment as research assistant (honorary). L. Liu thanks N. Crouseilles for the kind hospitality during her visit to University of Rennes. She acknowledges the support by National Key R\&D Program of China (2021YFA1001200), Ministry of Science and Technology in China, Early Career Scheme (24301021) and General Research Fund (14303022 \& 14301423) funded by Research Grants Council of Hong Kong.}}}
\author{Nicolas Crouseilles \thanks{Univ Rennes, Inria Rennes (Mingus team) and IRMAR UMR CNRS 6625, F-35000 Rennes, France \& ENS Rennes 
(\email{nicolas.crouseilles@inria.fr}).}, 
\and
Zhen Hao \thanks{School of Mathematics and Statistics, Wuhan University, Wuhan 430072, P. R. China (\email{zhhao\_math@whu.edu.cn}).}, 
\and
Liu Liu \thanks{The Chinese University of Hong Kong, Hong Kong 
(\email{lliu@math.cuhk.edu.hk}).}}

\headers{A bi-fidelity method for velocity space}{N. Crouseilles, Z. Hao, and L. Liu}

\makeatletter
\providecommand{\headerps@out}{}
\makeatother

\begin{document} 

\maketitle

\begin{abstract}
    In this paper, we introduce a bi-fidelity algorithm for velocity discretization of Boltzmann-type kinetic equations under multiple scales. The proposed method employs a simpler and computationally cheaper low-fidelity model to capture a small set of significant velocity points through the greedy approach, then evaluates the high-fidelity model only at these few velocity points and to reconstruct a bi-fidelity surrogate. This novel method integrates a simpler collision term of relaxation type in the low-fidelity model and an asymptotic-preserving scheme in the high-fidelity update step. Both linear Boltzmann under diffusive scaling and the nonlinear full Boltzmann in hyperbolic scaling are discussed. We show the weak asymptotic-preserving property and empirical error bound estimates. Extensive numerical experiments on linear semiconductor and nonlinear Boltzmann problems with smooth or discontinuous initial conditions and under various regimes have been carefully studied, which demonstrates the effectiveness and robustness of our proposed scheme. 
\end{abstract}



\input{1_Introduction}

\input{2_Kinetic_models}

\input{3_Bifidelity}

\input{4_Analysis}

\input{5_Numerics}

\input{6_Conclusion}

\input{7_Appendix}

\bibliography{ref}
\bibliographystyle{siamplain}

\end{document}

%% file: MMS_shared.tex
\usepackage{lipsum,mathtools}
\usepackage{amsfonts}
\usepackage{graphicx}
\usepackage{epstopdf}
\usepackage{algorithmic}
\ifpdf
  \DeclareGraphicsExtensions{.eps,.pdf,.png,.jpg}
\else
  \DeclareGraphicsExtensions{.eps}
\fi

\newsiamremark{remark}{Remark}
\newsiamremark{hypothesis}{Hypothesis}
\crefname{hypothesis}{Hypothesis}{Hypotheses}
\newsiamthm{claim}{Claim}
\headers{CEM-GMsFEM for Schr\"{o}dinger equations}{X. Jin, L. Liu, X. Zhong and E. T. Chung}
\title{Efficient numerical method for the Schr{\"o}dinger equation with high-contrast potentials\thanks{Submitted to the editors February 8, 2025.
\funding{L. Liu acknowledges the support by National Key R\&D Program of China (2021YFA1001200), Ministry of Science and Technology in China, Early Career Scheme (24301021) and General Research Fund (14303022 \& 14301423) funded by Research Grants Council of Hong Kong. E. Chung's research is partially supported by the Hong Kong RGC General Research Fund Projects 14305423 and 14305624.}}}
\author{Xingguang Jin\thanks{Department of Mathematics, The Chinese University of Hong Kong, Hong Kong
(\email{xgjin@math.cuhk.edu.hk}).}
\and Liu Liu\thanks{Department of Mathematics, The Chinese University of Hong Kong, Hong Kong
(\email{lliu@math.cuhk.edu.hk}).}
\and Xiang Zhong\thanks{Corresponding author. Department of Mathematics, The Chinese University of Hong Kong, Hong Kong
(\email{xzhong@math.cuhk.edu.hk}).}
\and Eric T. Chung\thanks{Department of Mathematics, The Chinese University of Hong Kong, Hong Kong
  (\email{tschung@math.cuhk.edu.hk}.})}

\usepackage{amsopn}


%% file: 1_Introduction.tex
\section{Introduction}

Bi-fidelity methods, or more generally multi-fidelity methods, leverage fast and inexpensive low-fidelity models together with accurate but costly high-fidelity models to solve complex problems~\cite{ZhuNarayanXiu2014, NarayanGittelsonXiu2014,ZhuXiu2017,Peherstorfer-Willcox-Gunzburger-2018}, thereby achieving reliable results at reduced expense.
These strategies have found particular success in uncertainty quantification, where they mitigate the burden of repeated high-cost simulations by exploiting correlations between model fidelities to guide efficient sampling and reconstruction. In recent years, bi-fidelity methods have proven effective in the field of uncertainty quantification (UQ), particularly when stochastic collocation (SC) methods are used.
The main challenge central to SC methods is the high simulation cost.
This is because, in many complex systems, an accurate high-fidelity simulation can
be computationally expensive. Furthermore, many stochastic algorithms such as SC require repetitive runs of the underlying deterministic simulation, which makes the overall stochastic simulation highly time consuming. The bi-fidelity method in SC seeks to leverage some approximate, low-fidelity models for the underlying problem that are computationally cheaper to run. In practice, the low-fidelity models usually contain simplified physics
\cite{NarayanGittelsonXiu2014} and/or are simulated on a much coarser physical mesh 
\cite{Liu-Zhu-2020}. By carefully combining results from both types of models, these methods can achieve high accuracy while greatly reducing computational cost, making them useful for problems where high-fidelity simulations alone would be too expensive.

Kinetic equations, which describe particle systems through distribution functions in phase space, exemplify such demanding problems. Their high-dimensional velocity domains, combined with intricate and stiff collision mechanisms, often lead to prohibitive computational loads, especially in multi-scale settings where classical solvers struggle with nonlinearity~\cite{DimarcoLiuPareschiZhu2021}. On the other hand, the presence of multiple scales and/or large velocities makes the kinetic equation stiﬀ. Classical stiﬀ solvers may be hard to use when we have to invert a very large nonlinear system. Traditional discretizations require dense velocity meshes to capture fine details, yet this amplifies costs without guaranteeing efficiency across regimes from kinetic to hydrodynamic limits.
To mitigate these computational burdens, various model reduction strategies have been proposed. Among the most prominent are hybrid solvers, which dynamically couple expensive kinetic solvers with computationally cheaper fluid models in different spatial domains~\cite{Tiwari-1998, Filbet-Rey-2015}. The core challenge for these hybrid methods lies in designing robust and efficient criteria for the decomposition of the spatial domain and managing the complex interfaces between the kinetic and fluid regions.

Our approach, the bi-fidelity or more general multi-fidelity method, is motivated by a similar goal of saving the computational cost but employs a different strategy. Instead of switching between different models based on domain decomposition and some criteria, we perform model reduction in the velocity space in the multi-fidelity framework.
Specifically, we develop a bi-fidelity method tailored to the velocity discretization of Boltzmann-type kinetic equations. By constructing a low-fidelity surrogate that captures the principal behavior in velocity space and selecting a small set of ``important'' velocity collocation points for high-fidelity evaluation, one can achieve a bi-fidelity approximation that balances cost and accuracy. At each time step, a low-fidelity model—typically a relaxation-based approximation—generates predictions across the full velocity grid, enabling a greedy algorithm to select a compact subset of pivotal velocity points. A high-fidelity model then computes accurate updates solely at these points, with the full solution reconstructed via projection coefficients derived from the low-fidelity data. This process ensures computational savings while preserving essential physical properties, such as convergence to macroscopic limits in stiff regimes through the asymptotic-preserving property of the bi-fidelity method.

\textbf{Main contribution:}
We apply this framework to two representative models: a linear semiconductor Boltzmann equation under diffusive scaling, featuring anisotropic collisions and external potentials, and a nonlinear Boltzmann equation with elastic particle interactions under a fluid scaling. Theoretical contributions include error estimates bounding the bi-fidelity deviation from high-fidelity references, proofs of weak asymptotic-preserving behavior, and an empirical metric for pre-assessing approximation quality without full high-fidelity runs. Numerical validations across various scenarios that range from smooth initial conditions to discontinuous shocks highlight the method's adaptability and robustness. Significantly, our approach demonstrates efficiency by \textit{adaptively} choosing representative velocities grid points(up to 1/25 the total number of grid points) in 1D2V nonlinear Boltzmann case while maintaining high accuracy, representing a substantial computational reduction. Moreover, this reduction also opens promising future work direction for extending the methodology to higher-dimensional problems (3D3V), which traditionally suffer from the curse of dimensionality and remain computationally prohibitive with conventional approaches.

The rest of the paper is organized as follows: Section 2 gives an introduction on the two Boltzmann type models. In Section 3, we discuss our designed bi-fidelity algorithm and give details of solving two model equations. In Section 4, we show the weak asymptotic-preserving property and an estimate on the empirical error bound. In Section 5, we conduct numerical experiments to demonstrate the efficiency and accuracy of our bi-fidelity approach. Lastly, we conclude the paper in Section 6.

%% file: 2_Kinetic_models.tex
\section{Boltzmann type kinetic models}

In kinetic theory of rarefied gases, 
the statistical description of a large particle system is given by a probability distribution function $f= f(t,x,v)$, characterizing the probability density of particles with velocity $v\in \mathbb{R}^{d_v}$, at position $x\in\Omega\subset\mathbb{R}^{d_x}$ and at time $t$.
The variation of $f$ is governed by a kinetic equation of the form
\begin{equation}
\label{vlasov_b}
    \partial_t f + v \cdot \nabla_x f + \nabla_x \Phi(t, x) \cdot \nabla_v f  =  \mathcal{Q}(f),
\end{equation}
where the left-hand side represents the advection of particles with an external electric potential $\Phi(t,x)$ acting on particles, and the right-hand side represents interactions between particles through a collision operator $\mathcal{Q}$, which acts only on the variable $v$.

In our numerical experiments, we will mainly study two models (called  Model 1 and Model 2) of the form \eqref{vlasov_b} that are presented below. 

\subsection{Model 1: semiconductor Boltzmann in diffusive scaling}
First, we consider a linear Boltzmann equation under the diffusive scaling in one dimensional space 
The distribution $f(t,x,v)$ satisfies the following kinetic equation~\cite{Markowich-Ringhofer-Schmeiser-1989}
 \begin{equation}
 \label{eq:diffusiveBP}
 \varepsilon \partial_t f + v \cdot \nabla_x f + \nabla_x \Phi \cdot \nabla_v f = \frac{1}{\varepsilon} \mathcal{Q}_{\text{LB}}(f), 
 \end{equation}
where $\varepsilon$ is the dimensionless mean free path  
and $\mathcal{Q}_{\text{LB}}$ is the anisotropic collision operator
\begin{equation}
    \mathcal{Q}_{\text{LB}}(f)(v) = \int_{\mathbb{R}^{d_v}} \sigma(v, w)\left(M(v)f(w) - M(w)f(v)\right) dw,
\label{equation:LinearBoltzmann}
\end{equation}
with $\sigma(v, w)$ being the collision cross-section and 
$\displaystyle M(v)= \frac{1}{\sqrt{2\pi}}
\exp(-|v|^2/2)$. This operator describes the collision of particles against a uniform background modeled by a Maxwellian distribution $M(v)$. Finally, we introduce the macroscopic density $\rho$ which is defined by $\displaystyle\rho = \int_{\mathbb{R}^{d_v}} f dv$. 

\subsection{Model 2: nonlinear Boltzmann equation}

Another model studied in this work is the nonlinear Boltzmann equation~\cite{Cercignani}. The model is given by
\begin{equation}
 \label{equation:NonlinearBoltzmann} \partial_t f + v\cdot\nabla_x f = \frac{1}{\e}\mathcal{Q}_{\text{NB}}(f, f), 
 \end{equation}
 The nonlinear operator $\mathcal{Q}_{\text{NB}}(f, f)$ describes elastic collisions between particles
\begin{equation*}
    \mathcal{Q}_{\text{NB}}(f, f)(v) = \int_{\mathbb{R}^{d_v}} \int_{\mathbb{S}^{d_v-1}} B\left(\left|v - v_*\right|, \sigma\right) \left(f(v') f(v_*') - f(v) f(v_*)\right)\mathrm{d}\sigma \, \mathrm{d}v_*, 
\end{equation*}
where $(v, v_*)$ and $(v', v_*')$ are the pre- and post-collision velocity pairs related by conservation of momentum and energy
\begin{equation*}
    v' = \frac{v + v_*}{2} + \frac{\left|v - v_*\right|}{2} \sigma, \quad
    v_*' = \frac{v + v_*}{2} - \frac{\left|v - v_*\right|}{2} \sigma,
\end{equation*}
with $\sigma$ the scattering direction varying on the unit sphere $\mathbb{S}^{d_v-1}$. The collision kernel $\displaystyle B\left(\omega;v, v_*, \sigma\right)=B\left(\omega;\left|v-v_*\right|, \cos \theta\right)$ is a non-negative function that only depends on $\omega$, $\left|v-v_*\right|$ and 
$\displaystyle\cos \theta=\frac{\sigma \cdot\left(v-v_*\right)}{\left|v-v_*\right|}$. The local Maxwellian $M(f)$ associated with $f$ is given by
\begin{equation}
    \label{equation:Maxwellian}
   M(f) = \frac{\rho}{2\pi T} \exp\left(-\frac{|v - u|^2}{2T}\right),
\end{equation}
where the density, bulk velocity and temperature are defined by
\begin{equation}
    \label{equation:macroscopic_moments}
\rho = \int_{\mathbb{R}^{d_v}} f \mathrm{d}v, \quad
u = \frac{1}{\rho}\int_{\mathbb{R}^{d_v}} v f\mathrm{d}v, \quad
T = \frac{1}{2\rho}\int_{\mathbb{R}^{d_v}} |v - u|^2 f \mathrm{d}v.
\end{equation}

%% file: 3_Bifidelity.tex
\section{The bi-fidelity method in velocity discretization of kinetic models}

In this section, we introduce the bi-fidelity method for the velocity discretization of collisional kinetic models like the ones presented above. For a review of the bi-fidelity method in UQ settings, see Appendix~\ref{appendix:UQ}.

\subsection{The general procedure}

We consider the probability density distribution \( f(t, x, v) \) on a velocity mesh $\mathcal{V}$, which is defined as a discrete set of points 
{\small
\begin{align*}
\mathcal{V} &= \Big\{ 
(v^1,\ldots,v^{d_v}) \in{\cal M}_{N_v, d_v}(\mathbb{R}) \; \Big| \; v^i = -L_v + k_i\Delta v, \\
& \;\;\;\;\;\;\;\; \hspace{4cm} 
k_i \in \{0,1,\ldots, N_v-1\}, \; i\in\{1,\ldots,d_v\} \Big\},
\end{align*}
}
where \(d_v\) represents the dimension of the velocity space, \(2L_v\) is the length of the velocity interval along each dimension, \(N_v\) is the number of velocity grid points in each dimension, \(\Delta v = 2L_v/N_v\) is the mesh spacing and $\mathcal{M}_{p,q}$ denotes the space of matrices of size $p \times q$ with real coefficients. 

\begin{figure}{}
    \centering
    \includegraphics[width=1\textwidth]{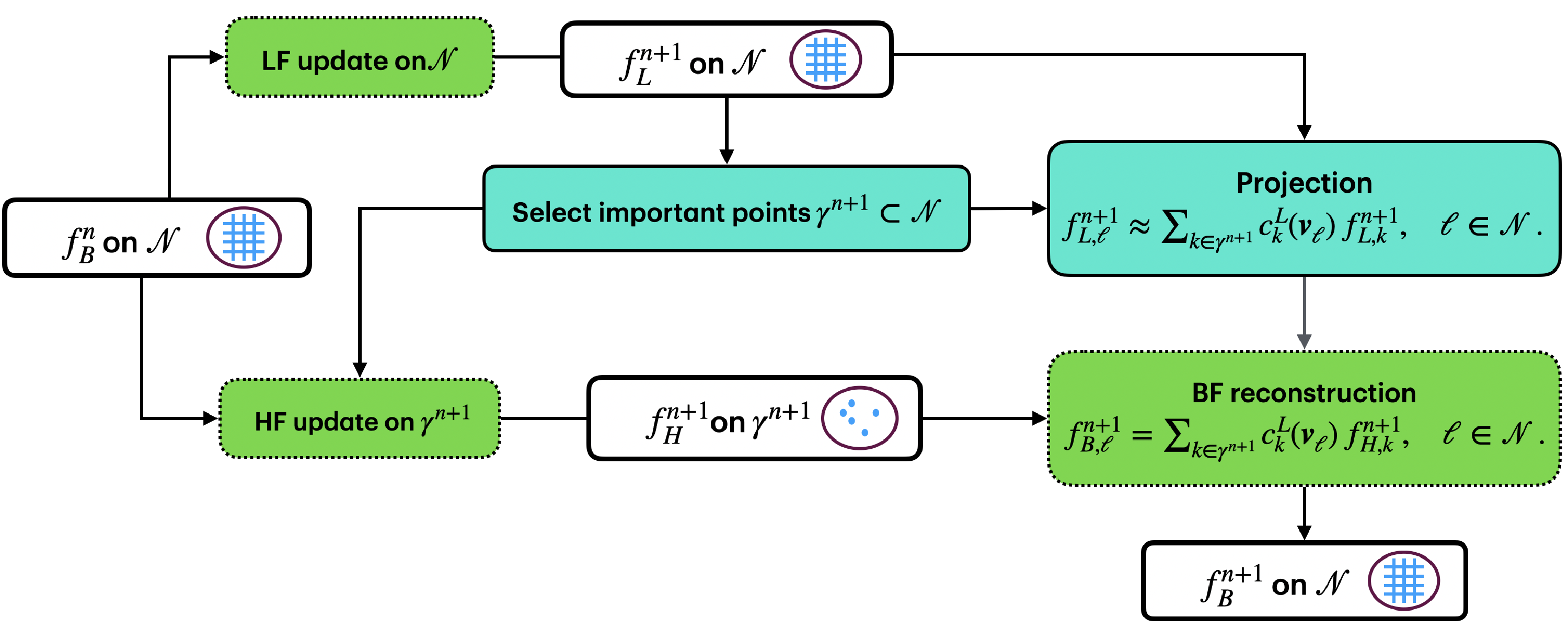}
    \caption{The bi-fidelity algorithm for the velocity discretization of kinetic models.}
    \label{fig:bifid}
\end{figure}

The bi-fidelity methodology centers around the idea of selecting `important' representative points out of a discrete set. In this regard, we index the $d_v$-dimensional mesh by a set $
    \mathcal{N} = \{ 1, \dots, N_v^{d_v} \} \text{ s.t } \mathcal{V} = \{ v_\ell \; | \; \ell\in\mathcal{N} \}$.
Note the difference between the subscript notation $v_\ell$ that is used to denote a velocity point in the mesh, and the superscript notation $v^i_\ell$ that is used to denote the $i$-th component of the velocity point $v_\ell$. We also define the cardinality of the velocity grid as $N = |\mathcal{N}| = N_v^{d_v}$.
For the spatial direction, we set the domain $[x_L, x_R]$ and define 
$
x_i = x_L + (i - 1) \Delta x, \; i\in \{1, ..., N_x\}$, 
where $\Delta x = (x_R - x_L)/N_x$ and $N_x$ the number of discretization points. 
At each time step $t^n$ and for a given position $x$, we treat the discrete solution vectors $\{f^n(x, v_\ell)\}_{\ell \in \mathcal{N}}$ as an ensemble of functions indexed by $\ell\in\mathcal{N}$. 
To obtain $\{f^{n+1}(x, v_\ell)\}_{\ell \in \mathcal{N}}$, the idea is to use a low-fidelity model to predict the solution first, then, based on this prediction, select a small subset of `important' velocity indices $\gamma^{n+1} \subset \mathcal{N}$ that forms a good basis for the entire set of velocity states, and finally use the high-fidelity model to predict the solution at the selected velocity indices.
The bi-fidelity method is summarized in Algorithm~\ref{alg:bifid} and Figure~\ref{fig:bifid}.

\textbf{Notations.} Throughout the paper, we will use $f^n_B$ to denote the numerical solution of a target equation at time $t^n$ obtained by the bi-fidelity method. The notations $f^{n+1}_L$ and $f^{n+1}_H$ are used to denote the low-fidelity and high-fidelity \textit{one-step predictions} from time $t^n$ to time $t^{n+1}$, respectively (see Algorithm~\ref{alg:bifid}). We also refer the notations $f^{n+1}_{B, \ell}$, $f^{n+1}_{L, \ell}$ and $f^{n+1}_{H, \ell}$ as numerical solutions at $v_\ell$, in particular $f^{n+1}_B(x, v_\ell)$, $f^{n+1}_L(x, v_\ell)$ and $f^{n+1}_H(x, v_\ell)$. For the convenience of notation, we let `BF', `HF' and `LF' stand for the bi-fidelity, high-fidelity and low-fidelity solutions, respectively. 

\begin{algorithm}
    \caption{A Bi-fidelity Time-Stepping Procedure}
    \label{alg:bifid}
    \begin{algorithmic}[1]
    \STATE \textbf{Input} Solution $\{f^n_{B, \ell}\}_{\ell \in \mathcal{N}}$ at time $t^n$.
    \STATE \textbf{Low-fidelity (LF) prediction:} Use a low-fidelity model to predict a solution $\{f_{L,\ell}^{n+1}\}$ for all velocity points $\ell\in\mathcal{N}$
    \STATE \textbf{Point selection:} Assemble the predicted LF solutions $\{f_{L,\ell}^{n+1}\}_{\ell \in \mathcal{N}}$ into a matrix
    {\small
    \begin{equation*}
    \begin{pmatrix}
        f_{L,1}^{n+1}(x_1) & \cdots & f_{L,N}^{n+1}(x_1)\\
        \cdot & \cdots & \cdot \\
        \cdot & \cdots & \cdot \\
        f_{L,1}^{n+1}(x_{N_x}) & \cdots & f_{L,N}^{n+1}(x_{N_x})
        \end{pmatrix}.
    \end{equation*}
    }Select a small subset of `important' velocity indices $\gamma^{n+1} \subset \mathcal{N}$ with $|\gamma^{n+1}|<\!\!\!< N$ using a greedy algorithm based on the pivoted Cholesky decomposition (see section~\ref{section:greedy}). 
    \STATE \textbf{Projection:} 
    Project the LF solution $\{f_{L,\ell}^{n+1}\}_{\ell \in \mathcal{N}}$ onto the basis spanned by the selected LF solutions $\{f_{L,k}^{n+1}\}_{k \in \gamma^{n+1}}$. In particular, find the coefficients $c^L_k(v_\ell)$ in
    $$ f_{L, \ell}^{n+1} \approx \sum_{k\in\gamma^{n+1}} c^L_{k}(v_\ell)f_{L,k}^{n+1},\quad \ell\in\mathcal{N},$$
    by solving the linear system
    \begin{equation}
        \label{equation:linear_system}
    \widetilde{\mathbf{G}}\mathbf{c} = \mathbf{f},
    \end{equation}
    where $\widetilde{\mathbf{G}} = \mathbf{G} - \mathbf{\Pi}$, $\mathbf{G}$ is the Grammian matrix of $(f^{n+1}_{L,k})_{k\in\gamma^{n+1}}$, 
    $\mathbf{c} = (c^L_{k}(v_\ell))_{k\in\gamma^{n+1},\; \ell\in\mathcal{N}}$ is the matrix of coefficients, and $\mathbf{f} = (\langle f_{L, k}^{n+1}, f_{L, \ell}^{n+1} \rangle)_{k\in\gamma^{n+1},\; \ell\in\mathcal{N}}$.
    $\mathbf{\Pi}$ denotes the projection matrix onto the kernel space of $\mathbf{G}$.
    \STATE \textbf{High-fidelity (HF) prediction:} Use the high-fidelity model to predict a solution $\{f_H^{n+1}(v_k)\}$ only for the selected `important' velocity points $k\in\gamma^{n+1}$.
    \STATE \textbf{Bi-fidelity (BF) reconstruction:} Use the LF projection coefficients but with the basis of the HF solutions to reconstruct the final bi-fidelity solution for all velocity points.
    $$ f_{B,\ell}^{n+1} = \sum_{k\in\gamma^{n+1}}c^L_{k}(v_\ell) f_{H,k}^{n+1}, \quad  \ell\in\mathcal{N}.$$
    \STATE \textbf{Output:} Solution $\{f^{n+1}_{B,\ell}\}_{\ell \in \mathcal{N}}$ at time $t^{n+1}$.
  \end{algorithmic}
\end{algorithm}

\begin{remark}
    In~\eqref{equation:linear_system}, the matrix $\widetilde{\mathbf{G}}$ is the Grammian matrix $\mathbf{G}$ minus the projection matrix $\mathbf{\Pi}$ onto the kernel space of $\mathbf{G}$. This is a common technique to avoid the inversion of a singular matrix.
\end{remark}

\subsection{Choice of low-fidelity and high-fidelity model}

The goal of the bi-fidelity method is to efficiently approximate a target problem by leveraging the accuracy of the high-fidelity model and the computational efficiency of the low-fidelity model. Hence, 
the best choice for the high-fidelity model used in Algorithm~\ref{alg:bifid} is simply a direct discretization of the target problem itself. 
To this end, we deploy asymptotic-preserving (AP) solvers for kinetic equations.
There have been many works developing robust numerical schemes for kinetic equations within the framework of asymptotic-preserving schemes; one can see a comprehensive review in~\cite{AP2022}. On the other hand, the low-fidelity model should be: $(i)$  a good approximation of the model problem and $(ii)$ computationally efficient. 

\subsection{Point selection: a greedy algorithm}
\label{section:greedy}

To form the set $\gamma^{n+1}$ of important velocity points, the idea is to search the velocity space using the low-fidelity data.
The searching process is usually conducted by a greedy algorithm.
In Algorithm~\ref{alg:cholesky}, we review the pivoted Cholesky greedy algorithm in~\cite{NarayanGittelsonXiu2014} for the selection of an index subset $\gamma\subset\mathcal{N}$ 
based on a given set of vectors $\mathbf{F} = \{f_\ell| \ell\in\mathcal{N}\}$. We set the stopping threshold $\delta$ which can be sufficiently small. 

\begin{algorithm}
\caption{Cholesky Decomposition Method for Selecting Interpolation Nodes from a Candidate Set $\mathcal{N}$.}
\label{alg:cholesky}
\begin{algorithmic}[1]
\STATE \textbf{Input:} ensemble $\mathbf{F} = \{f_1, f_2, \dots, f_{N}\}$, let $M$ the maximum number of points to be selected. 

\STATE Initialize index vector $\gamma = [1, \dots, M]$, inner product matrix $\mathbf{L} = \mathbf{0}_{N \times M}$. 

\STATE Initialize diagonal values $\mathbf{w}$ of the Gramian matrix of $\mathbf{F}$: $w_\ell = \|f_\ell\|_{\ell^2_x}^2$ for $\ell = 1, \dots, N$. 

\FOR{$n = 1, \dots, M$}
    \IF{$\max_{\ell \in \{n, \dots, N\}} w_\ell < \delta$}
        \STATE $n \gets n - 1$
        \STATE \textbf{break}
    \ENDIF
     \STATE $\gamma(n) \gets \arg\max_{\ell \in \{n, \dots, N\}} w_\ell$
    \STATE Exchange row $n$ and $\gamma(n)$ in $\mathbf{L}$, and exchange column $n$ and $\gamma(n)$ in $\mathbf{F}$ and in $\mathbf{w}$
     \STATE $r_\ell \gets \langle f_\ell, f_n \rangle_x - \sum_{j=1}^{n-1} L_{\ell,j} L_{n,j}$ for $\ell = n+1, \dots, N$
     \STATE $L_{n,n} \gets \sqrt{w_n}$
    \STATE $L_{\ell,n} \gets r_\ell / L_{n,n}$ for $\ell = n+1, \dots, N$
    \STATE $w_\ell \gets w_\ell - L_{\ell,n}^2$ for $\ell = n+1, \dots, N$
\ENDFOR
\STATE Truncate rows $M+1, \dots, N$ of $\mathbf{L}$
\STATE \textbf{Output:} selection of $\gamma$, Cholesky factor $\mathbf{L}$ for $\{ f_k | k\in\gamma\}$
\end{algorithmic}
\end{algorithm}

\subsection{Studied models}
\label{section:example_problems}

In this section, we consider two class of kinetic models that are numerically computed by our newly designed bi-fidelity scheme. 
In both problems, we propose to use the relaxation operator as a collision surrogate in the low-fidelity model: 
\begin{equation}
\label{equation:relaxation_operator}   \mathcal{Q}_{\text{LF}}:=
\lambda\left(\mathcal{E}_{[f]} - f\right), 
\end{equation}
where $\mathcal{E}_{[f]}$ is the local equilibrium state associated with the distribution function $f$. In the linear Boltzmann case, we let 
$\mathcal{E}_{[f]}=\rho M(v)$, whereas in the full Boltzmann case $\mathcal{E}_{[f]}=M(f)$ as defined in \eqref{equation:Maxwellian}.
The choice of $\lambda$ and discretized low-fidelity solver will be discussed in later section. 

\subsubsection{Semiconductor Boltzmann in diffusive scaling}

We consider the bi-fidelity method applied to the problem~\eqref{eq:diffusiveBP}.
\vspace{2mm}

\noindent\textbf{High-fidelity model.} 
The high-fidelity model we look for should be a direct discretization of~\eqref{eq:diffusiveBP}. 
We now define the following even and odd functions for $v>0$
\begin{equation}
    \label{eq:evenodd}
    \begin{aligned}
        r(t,x,v) &= \frac{f(t,x,v) + f(t,x,-v)}{2}, \quad j(t,x,v) = \frac{f(t,x,v) - f(t,x,-v)}{2\epsilon},
    \end{aligned}
\end{equation}
so that \eqref{eq:diffusiveBP} can be reformulated as 
\begin{equation}
\label{oddevendec}
    \left\{
    \begin{aligned}
        &\partial_t r + v \partial_x j + \partial_x\Phi\, \partial_v j = \frac{1}{\varepsilon^2} \mathcal{Q}_{\text{LB}}(r), \\
        &\partial_t j + \phi (  v \partial_x r + \partial_x\Phi\, \partial_v r ) = -\frac{1}{\varepsilon^2} \left[ \mu j  +  (1 - \varepsilon^2 \phi) (  v \partial_x r + \partial_x\Phi\, \partial_v r ) \right], 
    \end{aligned}
    \right.
\end{equation}
where $\mathcal{Q}_{\text{LB}}$ is given by \eqref{equation:LinearBoltzmann}  
and $\phi$ is a control parameter that guarantees the positivity of $r$ and $1 - \e^2 \phi$ and can be chosen as $\phi = \min(1, \frac{1}{\e^2})$~\cite{Jin-Pareschi-2000}.
The collision frequency $\mu(v)$ is given by
\begin{equation}\label{Lambda}
\mu(v) = \int_{\mathbb R} \sigma(v, w)M(w) dw. 
\end{equation}
The strategy of designing an AP scheme is to tackle the stiffness of the non-local collision operator by penalizing it with a relaxation operator
\begin{equation}
\label{equation:lowfid_operator_semi}
   P(g):=\lambda \left( \rho M(v) - g\right), \mbox{ with } \rho = \int_{\mathbb{R}} g dv, 
\end{equation}
where $\displaystyle\lambda = \max_v \mu(v)$. The AP scheme proposed in~\cite{Deng-2012} reads
\begin{equation}
    \label{scheme:highfid_semi}
    \left\{
    \begin{aligned}
    \frac{r^{n+1} - r^n}{\Delta t}  & +   v \partial_x j^n +  \partial_x \Phi \, \partial_v j^n 
    = \frac{1}{\varepsilon^2} \mathcal{Q}_{\text{LB}}(r^n) - \frac{1}{\varepsilon^2}P(r^n) + \frac{1}{\varepsilon^2}P(r^{n+1}), \\
    \frac{j^{n+1} - j^n}{\Delta t} &+ \phi ( v \partial_x r^n + \partial_x\Phi\, \partial_v r^n) \\ & = -\frac{1}{\varepsilon^2} \left( \mu j^{n+1}  +  (1 - \varepsilon^2 \phi) (v \partial_x r^{n+1} + \partial_x \Phi \,\partial_v r^{n+1})\right), 
    \end{aligned}
    \right.
\end{equation}
where 
$$
\mathcal{Q}_{\text{LB}}(f)(v_k) =  \Delta v  \sum_{\ell\in\mathcal{N}}\sigma(v_k, v_\ell)\left(M(v_k) f(v_\ell) - M(v_\ell) f(v_k)\right). 
$$
\vspace{2mm}

\noindent\textbf{Low-fidelity model.} 
Note that in the high-fidelity scheme~\eqref{scheme:highfid_semi}, the computationally intensive part is the evaluation of the nonlocal collision operator $\mathcal{Q}_{\text{LB}}(r^n)$, which only appears in the scheme for $r$. In the low-fidelity model, we use the relaxation operator $P$ from ~\eqref{equation:lowfid_operator_semi} to approximate the collision operator.
Hence, the bi-fidelity procedure is only needed to solve the scheme for $r$.
The low-fidelity scheme for $r$ reads
\begin{equation}
    \label{scheme:lowfid_semi}
    \begin{aligned}
        \frac{r^{n+1} - r^n}{\Delta t} + &  v \partial_x j^n + \partial_x\Phi \, \partial_v j^n = \frac{1}{\varepsilon^2} P(r^{n+1}). 
        \end{aligned}
\end{equation}
Now we adapt Algorithm~\ref{alg:bifid} to the even-odd framework and summarize our method in Algorithm~\ref{alg:bifid_semi}.

\begin{algorithm}
    \caption{Bi-fidelity Algorithm for the Semiconductor Boltzmann Equation}
    \label{alg:bifid_semi}
    \begin{algorithmic}[1]
    \STATE \textbf{Input:} Solution $\{f^n_{B,\ell}\}_{\ell \in \mathcal{N}}$ at time $t^n$.

    \STATE \textbf{Even-odd reformulation:} Get $r^n_B$ and $j^n_B$ by the even-odd decomposition ~\eqref{eq:evenodd}. 

    \STATE \textbf{LF prediction:} Use~\eqref{scheme:lowfid_semi} to predict $\{r_{L,\ell}^{n+1}\}$ for all velocity points $\ell\in\mathcal{N}$.

    \vspace{5pt}

    \STATE \textbf{Point selection:} Select $\gamma^{n+1} \subset \mathcal{N}$ using $r_{L,\ell}^{n+1}$ and Algorithm~\ref{alg:cholesky}. 

    \vspace{5pt}

    \STATE \textbf{Projection:} 
    Solve $c^L_{k}(v_\ell)$ via
    $ r_{L,\ell}^{n+1} = \sum_{k\in\gamma^{n+1}} c^L_{k}(v_\ell)r_{L,k}^{n+1},\quad \text{for } \ell\in\mathcal{N}$. 

    \vspace{5pt}

    \STATE \textbf{HF prediction:} Use the first equation in~\eqref{scheme:highfid_semi} to compute solutions 
    $\{r_{H,k}^{n+1}\}$ only for the selected `important' velocity points $k\in\gamma^{n+1}$.
    \vspace{5pt}

    \STATE \textbf{BF reconstruction:} Use the LF projection coefficients and HF basis solutions to reconstruct the BF solution for all velocity points: 
    \begin{equation}
        \label{rB-eqn}
     r_{B,\ell}^{n+1} = \sum_{k\in\gamma^{n+1}}c_k^L(v_\ell) r_{H,k}^{n+1}, \qquad \ell\in\mathcal{N}.  
    \end{equation}

   \STATE Update $j_{B,\ell}^{n+1}$ using the second equation in \eqref{scheme:highfid_semi} for $ \ell\in\mathcal{N}$. 

    \STATE \textbf{Solution output:}  
   $f^{n+1}_{B,\ell} = r^{n+1}_{B,\ell} + \varepsilon j_{B,\ell}^{n+1}$,\quad for $\ell\in\mathcal{N}$.
   \end{algorithmic}
\end{algorithm}

\subsubsection{Nonlinear Boltzmann equation}

In this part, we apply the bi-fidelity method to the nonlinear Boltzmann equation~\eqref{equation:NonlinearBoltzmann}.
\vspace{2mm}

\noindent\textbf{High-fidelity model.} 
For the high-fidelity model, we use the BGK penalty method \cite{Filbet-Jin-2010,Dimarco-Pareschi-2013} which is based on the following equivalent form of~\eqref{equation:NonlinearBoltzmann}
\begin{equation}
\label{equation:NonlinearBoltzmann_BGK}
    \partial_t f+v \cdot \nabla_x f=\frac{\mathcal{Q}_{\text{NB}}(f, f)-P(f)}{\varepsilon}+\frac{P(f)}{\varepsilon},
\end{equation}
where $P(f):=\lambda\left(M(f) -f\right)$ is the BGK collision operator \cite{BGK}. 

In this form, the first term on the right-hand side of \eqref{equation:NonlinearBoltzmann_BGK} is less stiff and is treated explicitly, while the second term is treated implicitly.
The high-fidelity scheme for $f$ reads
\begin{equation}
\label{equation:highfid_nonlinearBoltzmann}
    \frac{f^{n+1}-f^n}{\Delta t}+v \cdot \nabla_x f^n=
   \frac{\mathcal{Q}_{\text{NB}}(f^n, f^n)-P(f^n)}{\varepsilon}+\frac{P(f^{n+1})}{\varepsilon}.
\end{equation}
To obtain $M(f^{n+1})$ in $P(f^{n+1})$, 
we first multiply~\eqref{equation:highfid_nonlinearBoltzmann} by $m(v) = (1, v, |v|^2)^\intercal$ and integrate over the velocity space. The collision terms vanish due to conservation property of $\mathcal{Q}_{\text{B}}$ and $P(f)$. This leads to
\begin{equation}
    \left(\rho^{n+1}, \,\rho^{n+1}u^{n+1},\, \frac{1}{2} \rho^{n+1} (|u^{n+1}|^2 + 2 T^{n+1})\right)^\intercal = \int_{\mathbb{R}^{d_v}} m(v) \left(f^n - v\cdot \nabla_x f^n\right) \mathrm{d}v,
\end{equation}
where $\rho^{n+1}, u^{n+1}, T^{n+1}$ are the moments of $f^{n+1}$ defined in ~\eqref{equation:macroscopic_moments}, with $M(f^{n+1})$ given by ~\eqref{equation:Maxwellian}.
The parameter $\lambda$ is chosen as $\displaystyle\lambda = \max_v \mathcal{Q}_{\text{NB}}^-$ to preserve positivity~\cite{Filbet-Jin-2010, Yan-Jin-2012}, 
where
\begin{equation}
\mathcal{Q}^-_{\text{NB}}(f, f) = f(v)\int_{\mathbb{R}^{d_v}} \int_{\mathbb{S}^{d_v-1}} B(v - v_\star, \sigma) f(v_\star) \, \mathrm{d}\sigma \, \mathrm{d}v_\star\,.
\end{equation}
For the computation of the  nonlinear Boltzmann equation, most of the computational load is due to the multidimensional structure of the collisional integral $\mathcal{Q}_{\text{B}}$, which leads to an unpracticable computational cost with quadrature rules. To this end, we review the DVM method~\cite{Panferov-Heintz-2002} in Appendix~\ref{appendix:DVM}. This issue could be alleviated by the bi-fidelity approach as it reduces the query points of the collisional integral to just a few `important' ones.
\vspace{2mm}

\noindent\textbf{Low-fidelity model. } 
As mentioned in equation~\eqref{equation:relaxation_operator} 
we choose the BGK equation introduced by Bhatnagar, Gross and Krook~\cite{BGK} as our low-fidelity model, to alleviate the high computational cost of the Boltzmann collision operator. The low-fidelity scheme for $f$ reads
{\small
\begin{equation}
    \label{scheme:fL_nonlinearBoltzmann}
    \frac{f^{n+1}-f^n}{\Delta t}+v \cdot \nabla_x f^n=\frac{\lambda\left(M(f^{n+1})-f^{n+1}\right)}{\varepsilon}. 
\end{equation}
}In our numerical tests, we also choose
\begin{equation}
    \label{equation:lambda_nonlinearBoltzmann}
    \lambda = \max_v \mathcal{Q}_{\text{NB}}^{-}\,. 
\end{equation}
Now we summarize the numerical scheme for the nonlinear Boltzmann equation in  Algorithm~\ref{alg:bifid_nonlinearBoltzmann} below. 

\begin{algorithm}
    \caption{Bi-fidelity Algorithm for the Nonlinear Boltzmann Equation}    \label{alg:bifid_nonlinearBoltzmann}
    \begin{algorithmic}[1]
    \STATE \textbf{Input:} Solution $\{f^n_{B,\ell}\}_{\ell \in \mathcal{N}}$ at time $t^n$.

    \STATE \textbf{LF prediction:} Use \eqref{scheme:fL_nonlinearBoltzmann} to get $\{f_{L,\ell}^{n+1}\}$ for all velocity points $\ell\in\mathcal{N}$.
    \vspace{5pt}

    \STATE \textbf{Point selection:} Select $\gamma^{n+1} \subset \mathcal{N}$ using $f_{L,\ell}^{n+1}$ and Algorithm~\ref{alg:cholesky}. 

    \vspace{5pt}

    \STATE \textbf{Projection:} 
    Solve $c^L_{k}(v_\ell)$ via
    $f_{L,\ell}^{n+1} = \sum_{k\in\gamma^{n+1}} c^L_{k}(v_\ell)f_{L,k}^{n+1},\quad \text{for }\ell\in\mathcal{N}$. 

    \vspace{5pt}

    \STATE \textbf{HF prediction:} Use the scheme~\eqref{equation:highfid_nonlinearBoltzmann} to compute a solution $\{f_{H,k}^{n+1}\}$ only for the selected `important' velocity points $k\in\gamma^{n+1}$.

    \vspace{5pt}

    \STATE \textbf{Bi-fidelity (BF) reconstruction:} Use the LF projection coefficients and HF basis solutions to reconstruct the BF solution for all velocity points:
    \begin{equation}
        \label{fB-eqn}
      f_{B,\ell}^{n+1} = \sum_{k\in\gamma^{n+1}}c_k(v_\ell)f_{H,k}^{n+1}, \qquad \ell\in\mathcal{N}.
    \end{equation}

    \STATE \textbf{Solution output:} Gain $\{f^{n+1}_{B,\ell}\}_{\ell\in\mathcal{N}}$ at time $t^{n+1}$. 
 \end{algorithmic}
\end{algorithm}
\begin{remark}
In practice, to enforce the AP property of the high-fidelity prediction on a partial mesh, we adopt the BGK-penalty approach and incorporate it with a higher order discretization in time, namely Type A IMEX-RK modification~\cite{Dimarco-Pareschi-2013} in the high-fidelity update shown by
{\small
\begin{equation}
    \label{scheme:fH_nonlinearBoltzmann_typeA}
    \begin{aligned}
        f_{H,\ell}^{(1)} &= f_{B,\ell}^n + \frac{\Delta t}{\varepsilon} \lambda \left( M(f_{H}^{(1)})(v_\ell) - f_{H,\ell}^{(1)} \right), \quad \ell\in\mathcal{N}, \\[5pt]
        f_{H,k}^{(2)} &= f_{B,k}^n - \Delta t  v_k \cdot \nabla_x f_{H,k}^{(1)} + \frac{\Delta t}{\varepsilon} \lambda \left( M(f_{H}^{(2)})(v_k) - f_{H,k}^{(2)} \right) \\[5pt]
        &\quad + \frac{\Delta t}{\varepsilon} \Big[\mathcal{Q}_{\text{NB}}(f_H^{(1)}, f_H^{(1)})(v_k) - \lambda
        \left( M(f^{(1)}_{H})(v_k) - f_{H,k}^{(1)} \right)\Big], \quad k\in\gamma^{n+1}, \\[5pt]
        f_{H,k}^{n+1} &= f_{H,k}^{(2)}, \quad k\in\gamma^{n+1}.
    \end{aligned}
\end{equation}
}
\end{remark}

%% file: 4_Analysis.tex
\section{Error analysis and asymptotic-preserving property}
\label{sec:analysis}

Before proceeding to the numerical experiments, we analyze the error of the bi-fidelity reconstruction step. The following theorem provides a bound on the error between the bi-fidelity solution $f_B$ and the (uncomputed) full high-fidelity solution $f_H$. 
\begin{theorem}[Bi-Fidelity Approximation Error]
\label{thm:fB-fH-error}
Suppose we have the low- and high-fidelity predictions \( f^{n+1}_{L,\ell}(x_i) \), \( f^{n+1}_{H,\ell}(x_i) \), as well as the bi-fidelity reconstruction \( f^{n+1}_{B,\ell}(x_i) \) from the bi-fidelity algorithm for all \( \ell \in \mathcal{N} \) and all \( x_i \) with \( 1 \leq i \leq N_x \) at time \( t^{n+1} \). Suppose $\gamma^{n+1}$ is the set of selected nodes from the greedy algorithm at the current time step. Then, the error between the high-fidelity prediction and the bi-fidelity reconstruction at this step satisfies the following bound
$$
| f^{n+1}_{H,\ell}(x_i) - f^{n+1}_{B,\ell}(x_i) | \lesssim \max_{\ell \in \mathcal{N}} | f^{n+1}_{H,\ell}(x_i) - f^{n+1}_{L,\ell}(x_i) | + \delta, \quad \forall \ell \in \mathcal{N}, \quad \forall i = 1, \ldots, N_x,
$$
where \( \delta \) is the threshold in Algorithm~\ref{alg:cholesky}. In addition,
we have
$$
\| f^{n+1}_{H,\ell} - f^{n+1}_{B,\ell} \|_{\ell^2_{N_x}} \lesssim \sqrt{N_x}
\left( \max_{\ell \in \mathcal{N}} \| f^{n+1}_{H,\ell} - f^{n+1}_{L,\ell} 
\|_{\ell^2_{N_x}} + \delta \right), \quad \forall \ell \in \mathcal{N}. 
$$
Here we use the $\ell^2_N$ norm for a vector $u\in \mathbb{R}^N$: $\|u\|_{\ell^2_N}^2 = \sum_{i=1}^N u_i^2$.

\begin{proof}
For the convenience of notations, in the proof below we omit the supercripts of discrete solutions since they are all at time $t^{n+1}$.
First, similar to \cite{Gamba-Jin-Liu-2021}, we observe that 
{\small
\begin{equation}
\label{decomposition}
\begin{aligned}
    &f_{H,\ell}(x_i) - f_{B,\ell}(x_i) = f_{H,\ell}(x_i) - \sum_{k \in \gamma} c^L_k(v_\ell) f_{H,k}(x_i) \\
    = &f_{H,\ell}(x_i) - f_{L,\ell}(x_i) + \sum_{k \in \gamma} c^L_k(v_\ell) (f_{L,k}(x_i) - f_{H,k}(x_i)) + f_{L,\ell}(x_i) - \sum_{k \in \gamma} c^L_k(v_\ell) f_{L,k}(x_i) .
\end{aligned}
\end{equation}
}First, we establish an error bound on the last term. In previous works that studies UQ problems \cite{Cohen, Gamba-Jin-Liu-2021, Liu-Zhu-2020}, deriving a projection bound that involves an algebraic decay of $|\gamma|$
requires specific assumptions on the parameter space (random variable in their work), which in our current work changes to the velocity space, thus those assumptions may not be satisfied. 
Hence we derive a more practical bound based on linear algebra operations in Algorithm~\ref{alg:cholesky}. 

Consider Algorithm~\ref{alg:cholesky} at the \( (n_0 + 1) \)-th iteration, where \( \mathbf{V}_{n_0} \) denotes the subspace spanned by the \( n_0 \) vectors selected up to this iteration. The residual vectors, whose \( \ell^2_{N_x} \) norms are $w_\ell,\ \ell\in\{n_0+1, \dots, N\}$ in the algorithm, now belong to the set
\(
\mathcal{R}_{n_0} := \{ \tilde{f}_{L,\ell} := f_{L,\ell} - \Pi_{\mathbf{V}_{n_0}}[f_{L,\ell}] \mid \ell \in \mathcal{N}, \quad f_{L,\ell} \notin \mathbf{V}_{n_0} \},
\)
where \( \Pi_{\mathbf{V}_{n_0}} \) represents the orthogonal projection onto \( \mathbf{V}_{n_0} \) with respect to the \( \ell^2_{N_x} \) norm. The algorithm is designed to terminate when the following condition is satisfied: 
\begin{equation}
\label{Delta}
\max_{\tilde{f} \in \mathcal{R}_{n_0}} \|\tilde{f}\|_{\ell^2_{N_x}}  \leq \delta.
\end{equation}
Suppose the algorithm terminates at this step, then for \( f_{L,\ell} \in \mathbf{V}_{n_0} \) one has
\begin{equation}
\label{Projection_error}
f_{L,\ell} - \sum_{k \in \gamma} c^L_k(v_\ell) f_{L,k} \approx f_{L,\ell} - \Pi_{\mathbf{V}_{n_0}} [f_{L,\ell}] = 0.
\end{equation}
On the other hand, for \( f_{L,\ell} \notin \mathbf{V}_{n_0} \), we know that
$
\|f_{L,\ell} - \Pi_{\mathbf{V}_{n_0}}[f_{L,\ell}]\|_{\ell^2_{N_x}} \leq \delta, 
$
due to the definition of threshold $\delta$ given in \eqref{Delta}. 
Thus, for all \( \ell \in \mathcal{N} \) and \( i = 1, \ldots, N_x \),
\begin{equation}
\label{est2}
| f_{L,\ell}(x_i) - \sum_{k \in \gamma} c^L_k(v_\ell) f_{L,k}(x_i) | \leq \delta.
\end{equation}

For the second term in~\eqref{decomposition}, we need to derive a bound for \( c_k^L(v_\ell) \). Denote 
$\displaystyle \mathcal{P}_{\gamma} [f_{L,\ell}] = \sum_{k \in \gamma} c_k^L(v_\ell) f_{L,k}$.
For each \( \ell \in \mathcal{N} \),
\vspace{-10pt}
{\small
\begin{equation}
\label{PU_L1}
\begin{aligned}
\| \mathcal{P}_{\gamma} [f_{L,\ell}] \|_{\ell^2_{N_x}}^2 &= \sum_{m,n = 1}^{|\gamma|} c_m^L(v_\ell) c_n^L(v_\ell)\, 
\langle f_{L,m}, f_{L,n} \rangle_x = (\mathbf{c}^L)^{T} \widetilde{\mathbf{G}} \mathbf{c}^L
\geq \lambda_0 \|\mathbf{c}^L\|_{\ell^2_{|\gamma|}}^2\,,
\end{aligned}
\end{equation}
}where \( f_{L,\ell} = \{ f_{L,\ell}(x_i) \}_{i=1}^{N_x} \) and \( \widetilde{\mathbf{G}} = \mathbf{G} - \mathbf{\Pi} \) where \( \mathbf{G} \) is the Gramian matrix of 
$\{f_{L,k}(x_i)\}_{k \in \gamma}$, 
with $\mathbf{\Pi}$ the projection matrix onto the kernel space of $\mathbf{G}$. Here \( \lambda_0 > 0 \) is the minimum eigenvalue of $\widetilde{\mathbf{G}}$, and $\mathbf{c}^L$ is the vector composed of $c_m^L(v_\ell)$ for $m=1,\cdots, |\gamma|$, with $(\mathbf{c}^L)^{T} $ being its transpose.
It is obvious that
\begin{equation}\label{PU_L2}
\| \mathcal{P}_{\gamma} [f_{L,\ell}] \|_{\ell^2_{N_x}}^2 \leq \| f_{L,\ell} \|_{\ell^2_{N_x}}^2.
\end{equation}
Combining \eqref{PU_L1}--\eqref{PU_L2}, one gets
\(
\|\mathbf{c}^L(v_\ell)\|_{\ell^2_{|\gamma|}} \leq \frac{1}{\sqrt{\lambda_0}} \| f_{L,\ell} \|_{\ell^2_{N_x}}.
\)
Define $\displaystyle\|u\|_{\ell^\infty_{\mathcal{N}}} := \max_{\ell \in \mathcal{N}} |u_\ell|$ and $\displaystyle\|u\|_{\ell^\infty_{|\gamma|}} := \max_{\ell \in \gamma} |u_\ell|$. By Cauchy-Schwarz inequality, we have
\begin{equation}
\label{est1}
\begin{aligned}
& \quad |\sum_{k \in \gamma}  c^L_k(v_\ell) (f_{L,k}(x_i) - f_{H,k}(x_i)) | \\ &\leq \|\mathbf{c}^L(v_\ell)\|_{\ell_{|\gamma|}^2} \|f_L(x_i) - f_H(x_i)\|_{\ell_{|\gamma|}^2} \leq \frac{1}{\sqrt{\lambda_0}} \| f_{L,\ell} \|_{\ell^2_{N_x}} \|f_L(x_i) - f_H(x_i)\|_{\ell_{|\gamma|}^2} \\
&\leq \frac{|\gamma|^{1/2}}{\sqrt{\lambda_0}} \| f_{L,\ell} \|_{\ell^2_{N_x}}  \|f_{L}(x_i) - f_{H}(x_i)\|_{\ell^\infty_{|\gamma|}} \lesssim \|f_{L}(x_i) - f_{H}(x_i)\|_{\ell^\infty_\mathcal{N}}.
\end{aligned}
\end{equation}
Combine the bounds~\eqref{est1}, \eqref{est2} and the decomposition~\eqref{decomposition}, by the triangle inequality, we obtain
\(
|f_{H,\ell}(x_i) - f_{B,\ell}(x_i)| \lesssim \max_{\ell \in \mathcal{N}} |f_{L,\ell}(x_i) - f_{H,\ell}(x_i)| + \delta. 
\)
\end{proof}

\begin{remark}
We remark that in equation \eqref{Projection_error} ``$\approx$" holds due to practically the vectors $f_{L,k}$ may not be strictly linearly independent. Given the Cholesky decomposition in Algorithm \ref{alg:cholesky} is efficient enough, one can assume that the difference between the two projection terms on the left and right hand side of equation \eqref{Projection_error} is negligible.
From Theorem 4.1, we conclude that the error between bi-fidelity and high-fidelity solution is bounded by the error between low-fidelity and high-fidelity solution, plus a threshold used in Algorithm \ref{alg:cholesky} that can be set as small as we want. Since we have access to good low-fidelity model, this bi-fidelity error estimate is more than adequate.
\end{remark}

\end{theorem}

\subsection{AP property}

In this part, we will prove the weak AP property of the 
bi-fidelity method. The following theorem shows the weak AP property for the nonlinear Boltzmann equation \eqref{equation:NonlinearBoltzmann_BGK} but similar result can be obtained in the semiconductor Boltzmann case, which we put in Appendix C. 
\begin{theorem}
\label{theorem:AP}
Suppose $\delta < \varepsilon$
and let $(f_B^n)_n$ the numerical solution obtained from Algorithm 
\ref{alg:bifid_nonlinearBoltzmann}. Then, $(f^n_B)_n$ satisfies the following property: 
$$
\left|f^n_B - M(f^n_B) \right| = \mathcal{O}(\varepsilon) \Rightarrow
\left|f^{n+1}_B - M(f^{n+1}_B) \right| = \mathcal{O}(\varepsilon). 
$$
In particular, as $\varepsilon \rightarrow 0$, 
our numerical scheme automatically becomes a consistent discretization of the macroscopic limit equation, thus satisfies the AP property. 
\begin{proof}
First, we notice that
\begin{equation}
    \label{eq:fB-MB}
\begin{aligned}
    f^{n+1}_B - M(f^{n+1}_B) &= (f^{n+1}_B - f^{n+1}_H) + (f^{n+1}_H - M(f^{n+1}_B)) \\
    &= (f^{n+1}_B - f^{n+1}_H) + (f^{n+1}_H - M(f^{n+1}_H)) + (M(f^{n+1}_H) - M(f^{n+1}_B)).
\end{aligned}
\end{equation}
For the high-fidelity update, the first-order scheme \eqref{equation:highfid_nonlinearBoltzmann} is equivalent to
\begin{align}
\label{fH}
\begin{split}
    f^{n+1}_H&= \frac{\varepsilon}{\varepsilon + \lambda \Delta t} (f^n_B - \Delta t v \cdot \nabla_x f^n_B) + \frac{\lambda \Delta t}{\varepsilon + \lambda \Delta t} M(f^{n+1}_H) \\
    &\quad + \frac{\Delta t}{\varepsilon + \lambda \Delta t} \Big(\mathcal{Q}_{NB}(f^n_B, f^n_B) - \lambda (M(f^n_B) - f^n_B)\Big). 
\end{split}
\end{align}
The first term on the right-hand-side of \eqref{fH} is bounded by $\mathcal{O}(\varepsilon)$, and the 
second term can be written as $ M(f^{n+1}_H) - \frac{\varepsilon}{\varepsilon + \lambda \Delta t} M(f^{n+1}_H) =  M(f^{n+1}_H) + \mathcal{O}(\varepsilon)$. 

Regarding the third term, since by assumption $f^n_B - M(f^n_B) = \mathcal{O}(\varepsilon)$, we have $\mathcal{Q}_{NB}(f^n_B, f^n_B) = \mathcal{Q}_{NB}(M(f^n_B),M(f^n_B))+\mathcal{O}(\varepsilon)=\mathcal{O}(\varepsilon)$ so that
\begin{equation}
\label{q-p_B}
    \mathcal{Q}_{NB}(f^n_B, f^n_B) - \lambda (M(f^n_B) - f^n_B) = \mathcal{O}(\varepsilon).
\end{equation}
Thus, from \eqref{fH}, we deduce that
\begin{equation}
    \label{eq:fH-MH}
    f^{n+1}_H - M(f^{n+1}_H) = \mathcal{O}(\varepsilon).
\end{equation}
On the other hand, using Theorem~\ref{thm:fB-fH-error}, we have
\begin{equation}
    \label{eq:fB-fH-error}
    |f^{n+1}_H - f^{n+1}_B| \lesssim \max_v |f^{n+1}_H - f^{n+1}_L| + \delta.
\end{equation}
Since the Maxwellian function $M(f)$ is smooth, we also have, using \eqref{eq:fB-fH-error} 
\begin{equation}
    \label{eq:MH-MB-error}
    |M(f^{n+1}_H) - M(f^{n+1}_B)| \lesssim |f^{n+1}_H - f^{n+1}_B| \lesssim \max_v |f^{n+1}_H - f^{n+1}_L| + \delta.
\end{equation}
The low-fidelity scheme~\eqref{scheme:fL_nonlinearBoltzmann} can be rewritten as
\begin{equation}\label{fL}
    f^{n+1}_L= \frac{\varepsilon}{\varepsilon + \lambda \Delta t} (f^n_B - \Delta t v \cdot \nabla_x f^n_B) + \frac{\lambda \Delta t}{\varepsilon + \lambda \Delta t} M(f^{n+1}_L).
\end{equation} 
Subtracting the two schemes \eqref{fH} and \eqref{fL}, one gets
\begin{align}
\begin{split}
\label{Difference_fHL}
    f^{n+1}_H &- f^{n+1}_L =  \frac{\lambda \Delta t}{\varepsilon + \lambda \Delta t} M(f^{n+1}_H) + \frac{\Delta t}{\varepsilon + \lambda \Delta t} \Big(\mathcal{Q}_{NB}(f^n_B, f^n_B) - \lambda (M(f^n_B) - f^n_B) \Big) \\
    &\quad  - \frac{\lambda \Delta t}{\varepsilon + \lambda \Delta t} M(f^{n+1}_L)  
    = \frac{\Delta t}{\varepsilon + \lambda \Delta t} \Big(\mathcal{Q}_{NB}(f^n_B, f^n_B) - 
    \lambda (M(f^n_B) - f^n_B)\Big),  
\end{split}
\end{align}
since $M(f^{n+1}_L) = M(f^{n+1}_H)$. 
Using \eqref{q-p_B}, we get from \eqref{Difference_fHL}
\begin{equation}
    \label{eq:high_low_error}
    |f^{n+1}_H - f^{n+1}_L| = \mathcal{O}(\varepsilon).
\end{equation}
Finally, from~\eqref{eq:fB-MB} and using \eqref{eq:fH-MH}, \eqref{eq:fB-fH-error}, \eqref{eq:MH-MB-error}, we have
$
\left| f^{n+1}_B - M(f^{n+1}_B) \right| \lesssim |f_H^{n+1}-f_L^{n+1}| + \mathcal{O}(\varepsilon + \delta),  
$
which gives the result using~\eqref{eq:high_low_error} 
by choosing $\delta < \varepsilon$.
\end{proof}
\end{theorem}

\subsection{An empirical error bound}
\label{section:empirical_error_estimation}

In practical simulations, it is crucial to assess the quality of a chosen surrogate (LF) model and predict the error between HF and BF solutions before fully deploying the multi-fidelity algorithm. The theoretical error bound given in Theorem~\ref{thm:fB-fH-error} 
is often lack of practicality. To address this issue, the authors in \cite{Gao-Zhu-Wang-2020} introduced an empirical error estimation for UQ problems, which is easy to implement and effective in predicting the accuracy of bi-fidelity approximations without computing the HF solution in the entire parameter space. 

Inspired by this empirical error estimation strategy, we  will adapt it to our current study. The key idea is to use information from the LF model and the selected points to construct a {\it computable} error bound. This empirical bound provides a more realistic and practical measure of the predictive capability of BF models, and can be evaluated beforehand--that is, one does not require knowledge of the HF solution over the whole parameter space, which is most times not achievable in practice. In the following, we outline the methodology for our problem. 
\begin{lemma}
\label{lemma:empirical_error_bound}
    Suppose from time $t^n$ to $t^{n+1}$ we have selected the important points $\gamma^{n+1}$. Denote $U_H(\gamma^{n+1})$ as the subspace spanned by the HF solutions evaluated at $v\in\gamma^{n+1}$, and  $P_{U_H(\gamma^{n+1})}$ the projection operator onto $U_H(\gamma^{n+1})$, with the distance 
    $d^H(f^{n+1}_H(v), U_H(\gamma^{n+1}))=\|f^{n+1}_H(v) - P_{U_H(\gamma^{n+1})} f^{n+1}_H(v)\|_{\ell^2_{N_x}}$. 
Then the relative error between the one-step BF update $f^{n+1}_B$ and the (hypothetical) HF update $f^{n+1}_H$ on the whole velocity grid can be bounded as follows, for any $v \in \mathcal{V}$: 
{\small
\begin{align*}
&\frac{\|f^{n+1}_H(v) - f^{n+1}_B(v)\|_{\ell^2_{N_x}}}{\|f^{n+1}_H(v)\|_{\ell^2_{N_x}}}
\leq 
\underbrace{\frac{d^H(f^{n+1}_H(v), U_H(\gamma^{n+1}))}{\|f^{n+1}_H(v)\|_{\ell^2_{N_x}}}}_{\text{\textit{relative distance}}}
+ \underbrace{\frac{\|P_{U_H(\gamma^{n+1})} f^{n+1}_H(v) - f^{n+1}_B(v)\|_{\ell^2_{N_x}}}{\|f^{n+1}_H(v)\|_{\ell^2_{N_x}}}}_{\text{\textit{in-plane error}}} \\
&\quad 
= \max_v\left[\frac{d^H(f^{n+1}_H(v), U_H(\gamma^{n+1}))}{\|f^{n+1}_H(v)\|_{\ell^2_{N_x}}}\right]
\left(1 + \frac{\|P_{U_H(\gamma^{n+1})} f^{n+1}_H(v) - f^{n+1}_B(v)\|_{\ell^2_{N_x}}}{d^H(f^{n+1}_H(v), U_H(\gamma^{n+1}))} \right). 
\end{align*}
}
\end{lemma}

The proof is trivial and we refer readers to~\cite{Gao-Zhu-Wang-2020}. In the above error bound, note that we used the hypothetical HF update $f^{n+1}_H(v)$ on the entire velocity space, which is not necessary to compute, as shown in Algorithm 3.1. Hence, it is less useful and too expensive in practical simulation. To remove the dependency on HF values $f^{n+1}_H(v)$ on the whole velocity domain, we first define the {\it similarity ratio}
\begin{equation}
    \label{equation:rs}
    R_s = \max_v\left[\frac{d^H(f^{n+1}_H(v), U_H(\gamma^{n+1}))}{\|f^{n+1}_H(v)\|_{\ell^2_{N_x}}}\right] \left. \middle/ \max_v\left[\frac{d^L(f^{n+1}_L(v), U_L(\gamma^{n+1}))}{\|f^{n+1}_L(v)\|_{\ell^2_{N_x}}}\right].  \right.
\end{equation}

When $R_s \approx 1$, one can replace the above bound by the following simpler form
\begin{equation}
    \frac{\|f^{n+1}_H(v) - f^{n+1}_B(v)\|_{\ell^2_{N_x}}}{\|f^{n+1}_H(v)\|_{\ell^2_{N_x}}} \lesssim \max_v \left[\frac{d^L(f^{n+1}_L(v), U_L(\gamma^{n+1}))}{\|f^{n+1}_L(v)\|_{\ell^2_{N_x}}} \right] \left(1 + R_e(v) \right),
\end{equation}
where the quantity $R_e(v)$ is defined by
\begin{equation}
    \label{equation:re}
    R_e(v) = \frac{\|P_{U_H(\gamma^{n+1})} f^{n+1}_H(v) - f^{n+1}_B(v)\|_{\ell^2_{N_x}}}{\|d^H(f^{n+1}_H(v), U_H(\gamma^{n+1}))\|_{\ell^2_{N_x}}},
\end{equation}
which measures the balance between the in-plane error and the relative distance. In addition, similar to~\cite{Gao-Zhu-Wang-2020}, we can use the last selected velocity point $v_{\text{end}} \in \gamma^{n+1}$ to serve as the testing point in this error surrogate for the BF approximation in the entire velocity space. We summarize this empirical error bound as follows: 

\vspace{5pt}

\noindent\textbf{Empirical error bound.}
\textit{
    If the similarity ratio $R_s \approx 1$, the empirical error bound which is calculated by
    \begin{equation}
        \label{equation:empirical_error}
        \mathcal{E}^{n+1} := \max_v \left[\frac{d^L(f^{n+1}_L(v), U_L(\gamma^{n+1}))}{\|f^{n+1}_L(v)\|_{\ell^2_{N_x}}} \right] \left(1 + R_e(v_{\text{end}}) \right)
    \end{equation}
    will serve a reliable bound of the true error, since
\begin{equation}
\label{equation:empirical_error_bound}
   \text{True Error}:=\frac{\|f^{n+1}_H(v) - f^{n+1}_B(v)\|_{\ell^2_{N_x}}}{\|f^{n+1}_H(v)\|_{\ell^2_{N_x}}} \lesssim \mathcal{E}^{n+1}. 
\end{equation}
}

From~\eqref{equation:empirical_error_bound}, we see that compared to the ``True Error" which requires computation of high-fidelity updated solution in the whole velocity space, the empirical bound $\mathcal{E}^{n+1}$ can be computed using only a few  selected velocity points, which is quite efficient. 
It can be seen that besides $R_s \approx 1$ (LF and HF models should be similar), the approximation quality of the BF approximation also depends on $R_e$, which describes the balance between the in-plane error and projective distance. A large $R_e$ indicates that the in-plane error is dominant over the projective distance leading to not good BF approximation. 
In summary, the effectiveness of the empirical error bound relies on two main ingredients: 
(i) when the similarity ratio $R_s$ is close to $1$, the empirical error bound given in~\eqref{equation:empirical_error_bound} serves as a reliable approximation of the true error. (ii) provided that $R_e$ remains moderately small, 
the bi-fidelity method should be capable of delivering satisfactory results.

%% file: 5_Numerics.tex
\section{Numerical Examples}

In this section, we present numerical examples to demonstrate the effectiveness and robustness of our designed bi-fidelity method. We consider numerical example based on two models mentioned in section~\ref{section:example_problems}: the nonlinear Boltzmann equation in hyperbolic scaling \eqref{equation:NonlinearBoltzmann} and the semiconductor Boltzmann equation in diffusive scaling \eqref{eq:diffusiveBP}. In the points selection step shown in Algorithm 3.2, the tolerance is set to be $\delta = 10^{-12}$.

\subsection{Smooth problems}

In this section, we conduct numerical experiments for problems with smooth initial conditions (I.C.). 

\subsubsection{Test I (a): Smooth, equilibrium I.C. for Model 1}

\begin{figure}
    \centering
    \begin{subfigure}[b]{0.3\textwidth}
        \includegraphics[width=\textwidth]{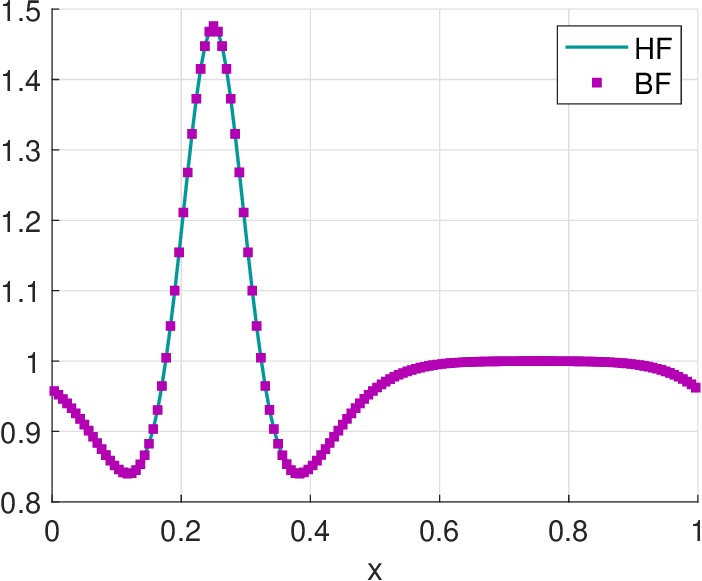}
        \caption{$\varepsilon = 1$}
        \label{fig:smooth_eps1}
    \end{subfigure}
    \hspace{2em}
    \begin{subfigure}[b]{0.3\textwidth}
        \includegraphics[width=\textwidth]{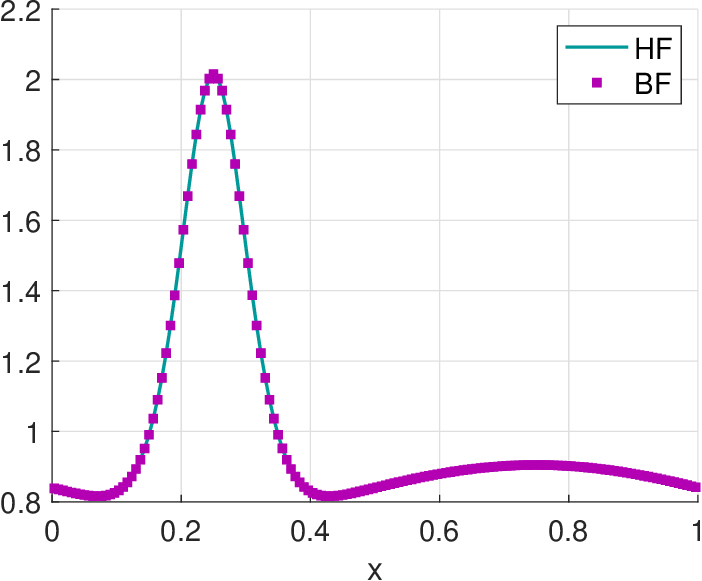}
        \caption{$\varepsilon = 10^{-1}$}
        \label{fig:smooth_eps1e-1}
    \end{subfigure}

    \vspace{2pt}

    \begin{subfigure}[b]{0.3\textwidth}
        \includegraphics[width=\textwidth]{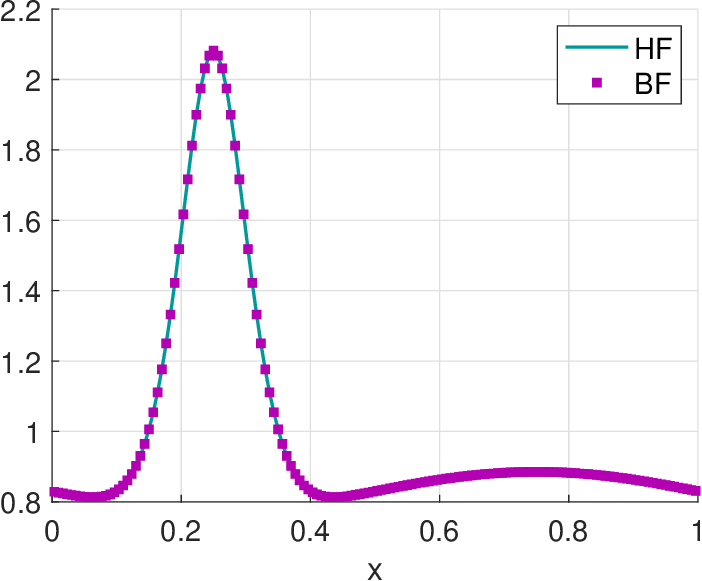}
        \caption{$\varepsilon = 10^{-3}$}
        \label{fig:smooth_eps1e-4}
    \end{subfigure}
    \hspace{2em}
    \begin{subfigure}[b]{0.3\textwidth}
        \includegraphics[width=\textwidth]{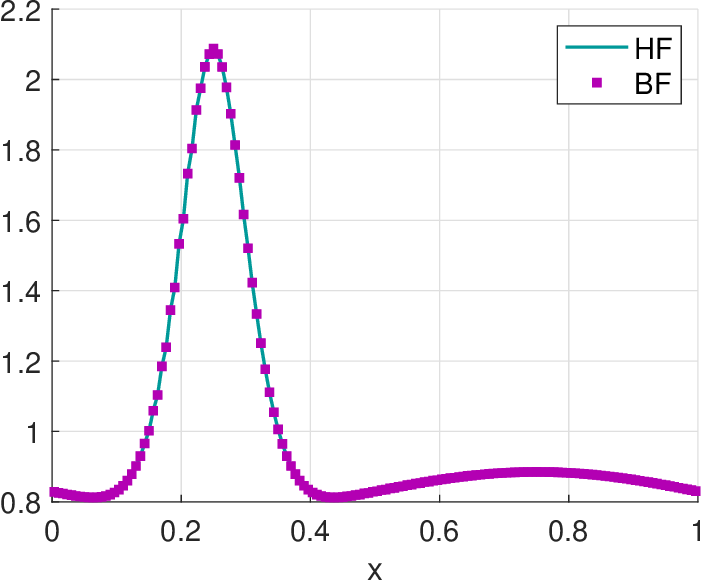}
        \caption{$\varepsilon = 10^{-8}$}
        \label{fig:smooth_eps1e-9}
    \end{subfigure}
    \caption{Test I (a): Density profiles computed from the bi-fidelity (dotted lines) and the high-fidelity (solid lines) reference solution for varying $\varepsilon$. }
    \label{fig:BPdensity}
\end{figure}

In this problem, we consider~\eqref{eq:diffusiveBP}
with external potential 
\( \Phi(x) = e^{-50e(1/4-x)^2} \) and initial condition 
$$
f^{\text{in}}(x, v) = \frac{1}{\sqrt{2\pi}}e^{\frac{-v^2}{2}}.
$$
The spatial domain is set to be $x\in [0,1]$ and periodic boundary condition is assumed. For collision cross section, we let
$\sigma(v, w) = 1 + M(v)e^{-(|v|^2 - |w|^2 + 1)^2} + M(w)e^{-(|v|^2 - |w|^2 - 1)^2}$.
In this test, we use $N_x = 150$, $N_v = 100$, $\Delta t = 5\times 10^{-5}$. First, we investigate the performance of BF approximation by comparing it with the HF reference solution. Here we distinguish between the `HF reference solution' and the `HF update' in the BF method, where the former is computed by simulating the HF solver throughout the entire computational time, while the latter is simply the intermediate solution updated in our BF Algorithm~\ref{alg:bifid}. The density profiles $x\mapsto \rho(t=0.1, x)$ for varying $\varepsilon$ are plotted in 
Figure~\ref{fig:BPdensity}. We conclude that with $\varepsilon \in \{ 1, 10^{-1}, 10^{-3}, 10^{-8}\}$, the BF approximation is consistent with the HF solution at each spatial point. The BF method produces accurate solution at the final time. 

\begin{figure}
    \centering
    \begin{subfigure}[b]{0.4\textwidth}
        \includegraphics[width=\textwidth]{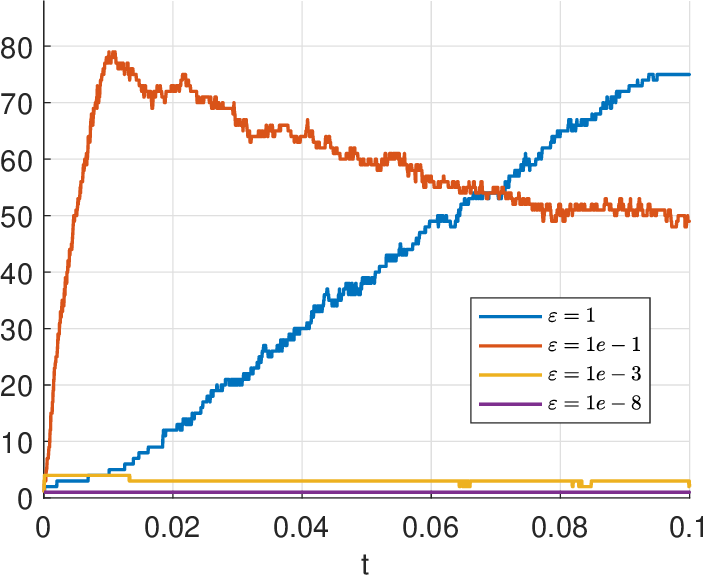}
    \end{subfigure}
    \hfill
    \begin{subfigure}[b]{0.4\textwidth}
        \includegraphics[width=\textwidth]{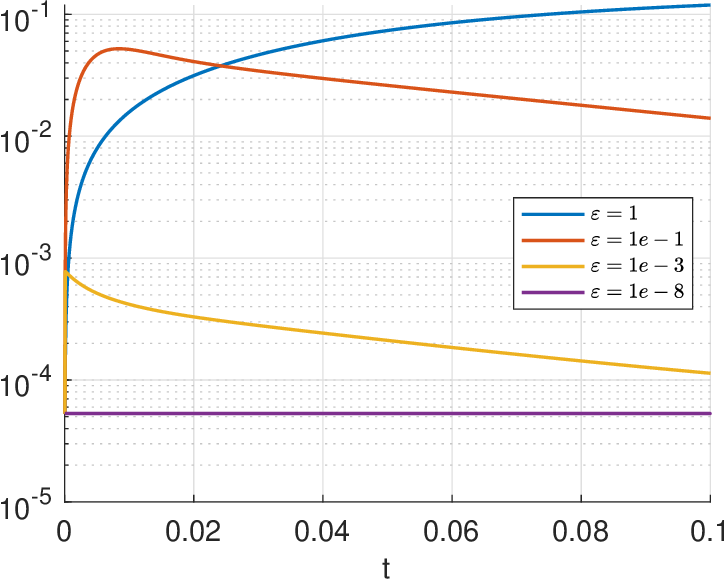}
    \end{subfigure}
    \caption{Test I (a): Left: Time evolution of $|\gamma|$ for different $\varepsilon$. Right: Time evolution of $\|f_B - \rho(f_B)M\|_{\ell^1}$ for different $\varepsilon$.}
    \label{fig:BPhistory}
\end{figure}

In Figure~\ref{fig:BPhistory}, we plot the time evolution of the size of important velocity point set $|\gamma|$ in the bi-fidelity algorithm, together with the time evolution of the distance between bi-fidelity solution and the local equilibrium, i.e., $\|f_B - \rho(f_B)M\|_{\ell^1}$.
Based on the two figures, we observe that as time propagates our algorithm is able to \textit{adaptively} select more points in the velocity space when $\varepsilon$ is larger, since in the kinetic regime, the solution is farther from equilibrium within the initial period of time and more velocity points are required to represent it. Under scenarios where $\varepsilon$ is small and/or the solution approaches the equilibrium state, fewer velocity points are needed. This phenomenon matches with our intuition since an increasing complexity of the bi-fidelity approximation in the velocity space is expected to bring about more significant points chosen in the algorithm. 

\begin{figure}
    \centering
    \includegraphics[width=0.4\textwidth]{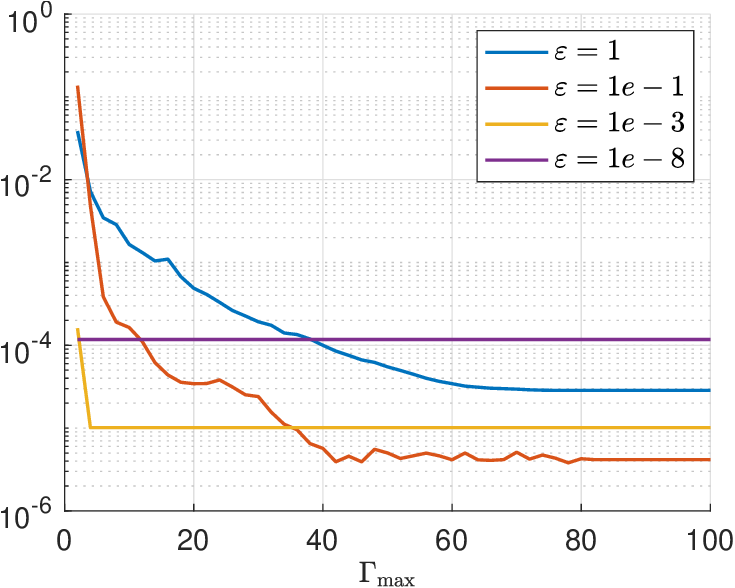}
    \caption{Test I (a): Relative error $\mathscr{E}$ at $t=0.1$ with respect to $\Gamma_{\text{max}}$. }
    \label{fig:BP_maxSelect_sweep}
\end{figure}

On the other hand, we study the convergence behavior of the BF approximation to the HF reference solution as the maximum allowed velocity points increases. In other words, we limit the size of the set $\gamma$ by a number $\Gamma_{\max}$,  
gradually increase $\Gamma_{\max}$ to observe the convergence between BF and HF solutions. 
In Figure~\ref{fig:BP_maxSelect_sweep}, the relative error in $\ell^1$ norm at final time $t=0.1$ against the above mentioned $\Gamma_{\max}$ is shown.   
The relative error in $\ell^1$ norm of phase space is defined by
{\small
\begin{equation}
\label{equation:error_vs_maxSelect}
    \mathscr{E}(t) = \frac{\left\|f_H(t,\cdot) - f_B(t,\cdot)\right\|_{\ell^1}}{\left\|f_H(t,\cdot)\right\|_{\ell^1}}.
\end{equation}
}We observe the relative error decreases as $\Gamma_{\max}$ gets larger and eventually saturates--adding more velocity points will not further increase the accuracy of BF solution, by then the error is dominated by other numerical discretizations in time and space. 

Moreover, when the errors saturate at the level of $\mathcal{O}(10^{-4})$ to $\mathcal{O}(10^{-6})$, the smaller $\varepsilon$ is, the fewer number of points are needed. This is due to the fact that when $\varepsilon$ is very small the model converges to the limiting diffusion equation and the distribution is close to its local equilibrium, thus less velocity points are required to describe the solution behavior in the velocity space. In a word, Figure~\ref{fig:BP_maxSelect_sweep} indicates that the greedy algorithm can effectively select the important velocity points for semiconductor Boltzmann equation varying from kinetic to diffusive regime. 

\subsubsection{Test I (b): Smooth, non-equilibrium I.C. for Model 2}

In this test, 
we consider the nonlinear Boltzmann equation~\eqref{equation:NonlinearBoltzmann} with two-dimensional velocity variable and one-dimensional spatial space. 

The computation of $\mathcal{Q}_{\text{B}}$ is performed with the DVM method~\cite{Panferov-Heintz-2002}, with details reviewed in Appendix~\ref{appendix:DVM}. The computational domain for velocity is set as \( v \in [-L_v, L_v]^2 \) with $L_v = 8$, we let $N_v = 32$ points in each dimension. The spatial domain is \( x \in [0, 1] \), we choose $N_x = 50$ grid points and assume periodic boundary conditions. The time step is chosen to satisfy the transport CFL condition $\Delta t = \frac{\Delta x}{2 L_v}$. 
We consider the Maxwell molecule model where the collision kernel $B$ is a constant, with $B = \frac{1}{2\pi}$.
We let the final time $t=0.2$ and show results for different Knudsen number with $\e \in \{ 1, 10^{-4}, 10^{-8}\}$ spreading from kinetic to hydrodynamic regimes.

First, we consider a smooth yet far-from-equilibrium initial condition which is composed of the sum of two Maxwellian distributions with different mean velocities:
\begin{equation*}
    f^{\text{in}}(x, v) = \frac{\rho^0(x)}{4\pi T^0(x)} \exp \left(- \frac{(|v - u^0(x)|)^2}{2T^0(x)}\right) + \frac{\rho^0(x)}{4\pi T^0(x)} \exp \left(- \frac{(|v + u^0(x)|)^2}{2T^0(x)}\right), 
\end{equation*}
where the macroscopic profiles are given by
$\rho^0(x) = \frac{2 + \sin (2\pi x)}{3}$, $u^0(x) = (\cos (2\pi x), 0)$, $T^0(x) = \frac{3 + \cos (2\pi x)}{4}$.

\begin{figure}
    \centering
    \begin{subfigure}[b]{0.8\textwidth}
        \caption{$\varepsilon = 1$}   \includegraphics[width=\textwidth]{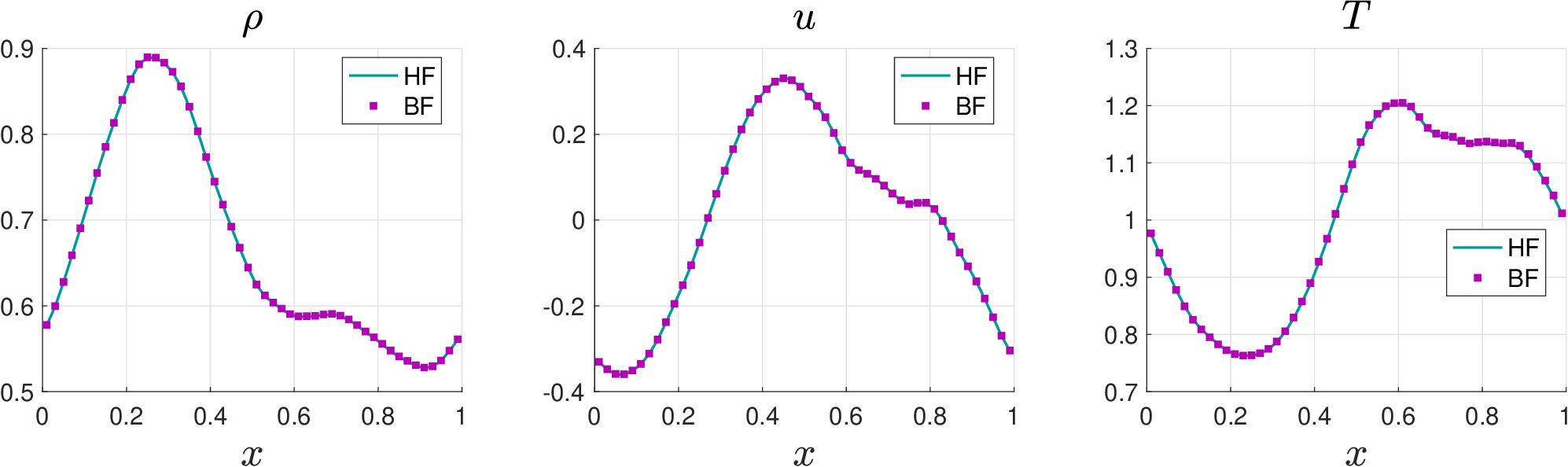}
    \end{subfigure}
    \hfill
    \hspace{1em}
    \begin{subfigure}[b]{0.8\textwidth}
        \caption{$\varepsilon = 10^{-4}$}        \includegraphics[width=\textwidth]{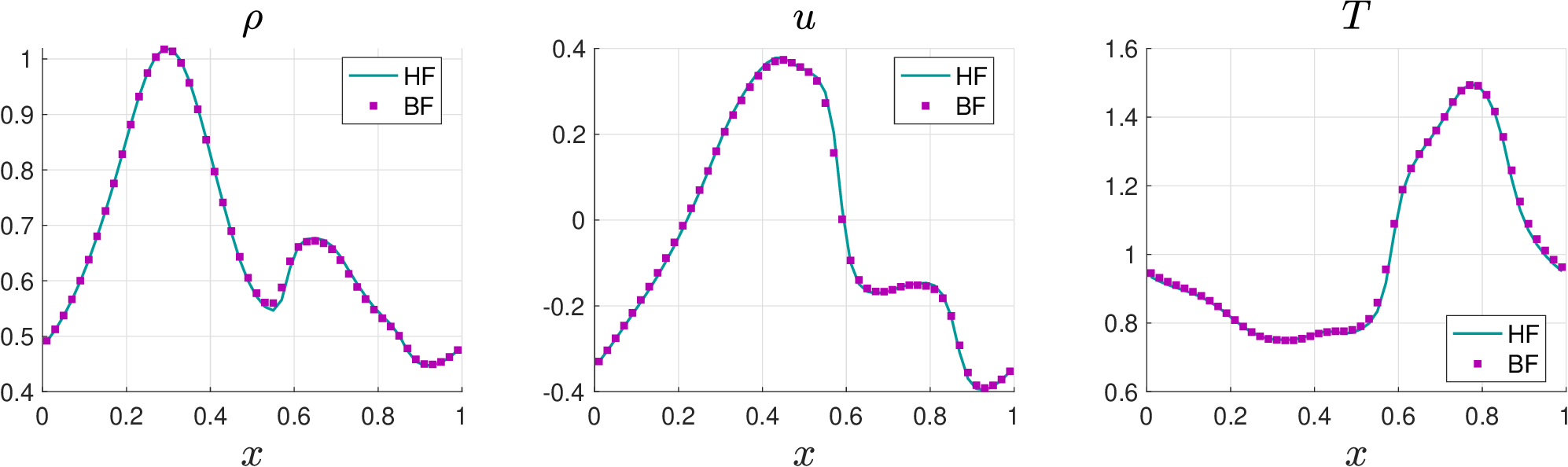}
    \end{subfigure}
    \hfill
    \hspace{1em}
    \begin{subfigure}[b]{0.8\textwidth}
    \caption{$\varepsilon = 10^{-8}$}        \includegraphics[width=\textwidth]{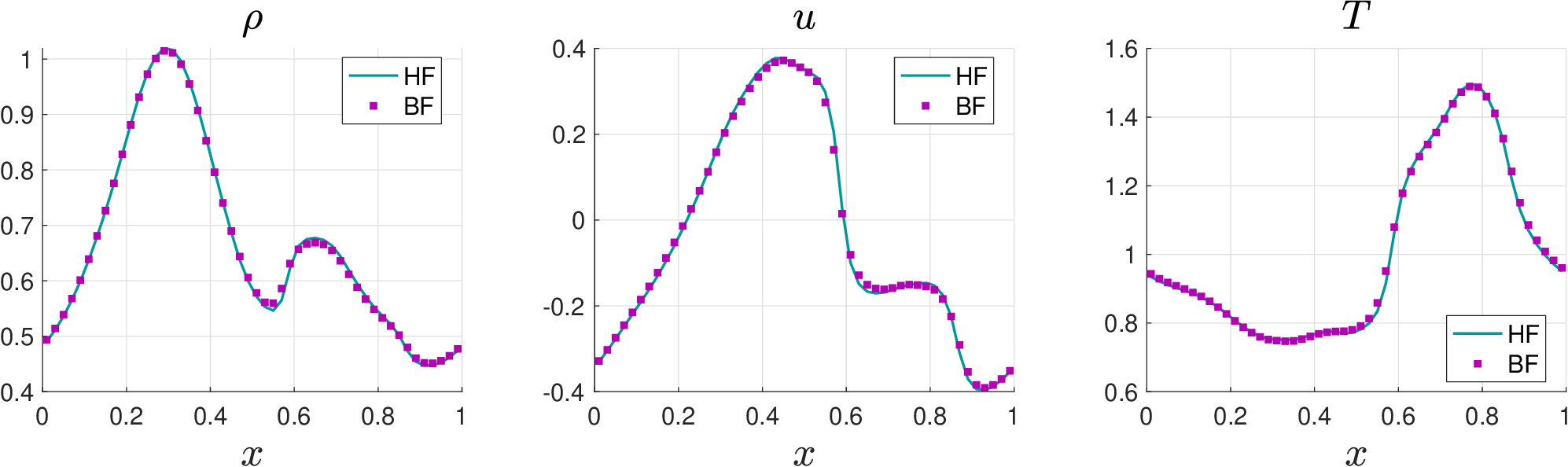}
    \end{subfigure}
    \caption{Test I (b): Macroscopic quantity profiles at $t = 0.2$ of bi-fidelity (dotted lines) and high-fidelity (solid lines) solutions for varying Knudsen numbers $\varepsilon=1$ (top row), $\varepsilon=10^{-4}$ (middle row) and $\varepsilon=10^{-8}$ (bottom row). Left: density $x\mapsto\rho(t=0.2, x)$. Middle: bulk velocity  $x\mapsto u(t=0.2, x)$. Right: temperature   $x\mapsto T(t=0.2, x)$.} 
    \label{fig:Boltz_smooth_data}
\end{figure}

\begin{figure}
    \centering
    \begin{subfigure}[b]{0.4\textwidth}        \includegraphics[width=\textwidth]{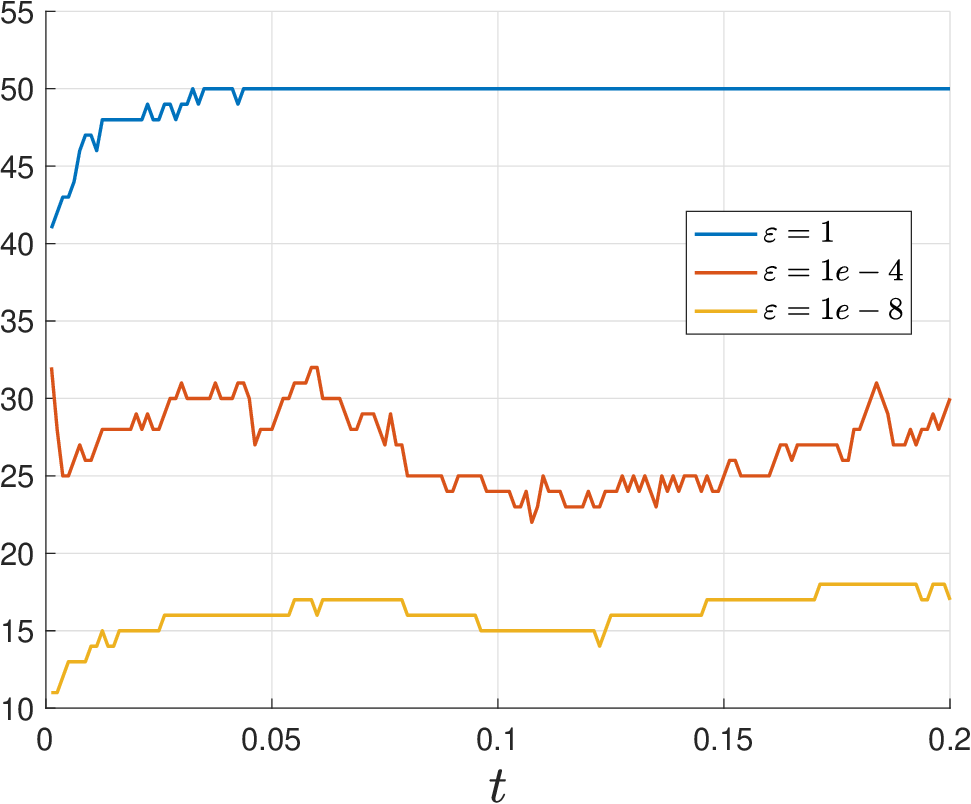}
    \end{subfigure}
    \hfill
    \begin{subfigure}[b]{0.4\textwidth}        \includegraphics[width=\textwidth]{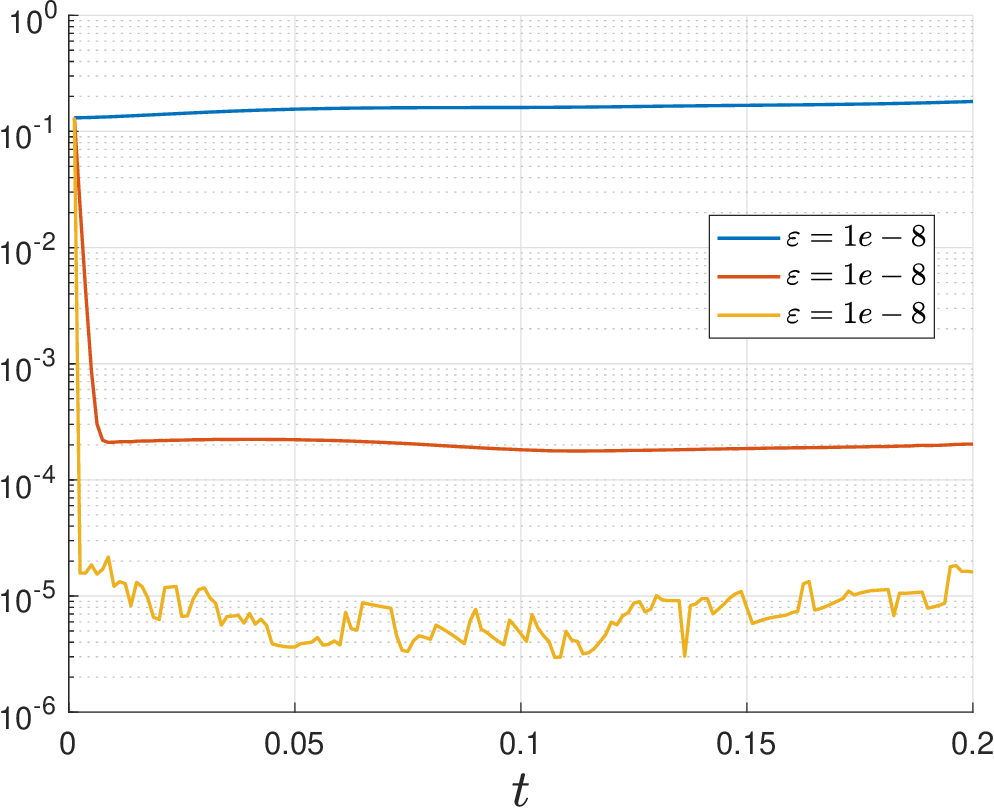}
    \end{subfigure}
    \caption{Test I (b): Left: Time evolution of $|\gamma|$ for different $\varepsilon$. Right: Time evolution of $\|f_B - M(f_B)\|_{\ell^1}$ for different $\varepsilon$. }    \label{fig:Boltz_smooth_history}
\end{figure}

In Figure~\ref{fig:Boltz_smooth_data}, it is shown that the macroscopic quantities (density $\rho$, bulk velocity $u$ and temperature $T$) computed from the BF solution are in excellent agreement with the HF reference solution at each spatial point and across different regimes. The time  evolution of the number of selected velocity points $|\gamma|$ and the distance from $f_B$ to equilibrium 
$\|f_B - M(f_B)\|_{\ell^1}$ is shown in Figure~\ref{fig:Boltz_smooth_history}. Similar to result in Test I (a), the BF algorithm is able to adaptively select more points for larger $\varepsilon$ when distribution solution is far  from its local equilibrium, whereas less points are picked when the solution approaches a hydrodynamic state. 
We note that among the whole velocity grid with size $N_v^2=1024$, the number of selected points is less than $50$, which significantly reduces the computational cost while capturing the solution behavior well. 

\begin{figure}
    \centering
\includegraphics[width=0.4\textwidth]{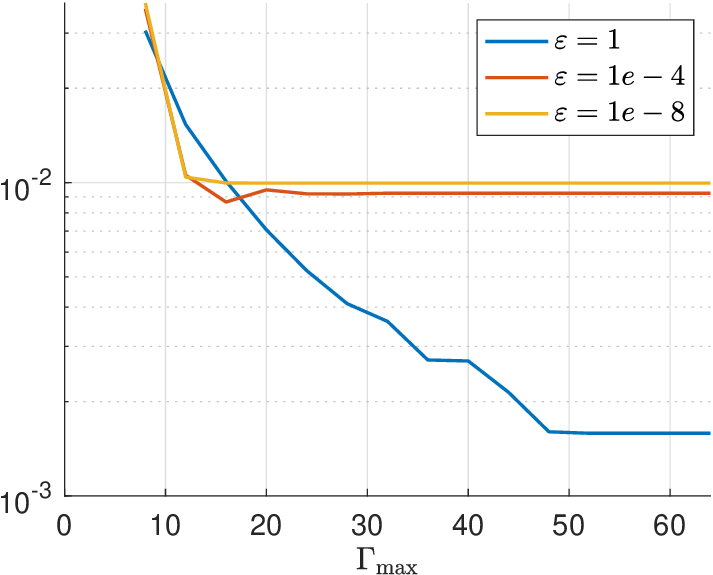}
\caption{Test I (b): Relative $\mathscr{E}(t=0.2)$ error with respect to $\Gamma_{\text{max}}$. }
    \label{fig:Boltz_smooth_error}
\end{figure}

In Figure~\ref{fig:Boltz_smooth_error}, we investigate the approximation error~\eqref{equation:error_vs_maxSelect} as a function of the maximum allowed set size $\Gamma_{\text{max}}$. Similar to Test I (a), one can see that the error decreases and eventually saturates as more points are included. For small $\varepsilon$, the distribution is close to equilibrium and the relative error saturates at less than $|\gamma|=20$ points whereas when 
 $\varepsilon = 1$, the relative error saturates until more than $|\gamma|=50$ points. 

\subsection{Discontinuous problems}

In this section, the bi-fidelity method is used to solve more challenging problems, with discontinuous initial conditions. The two benchmark tests are Riemann problem and blast wave problem \cite{Filbet-Jin-2010}. 

\subsubsection{Test II: Riemann problem for Model 2}

We consider a Riemann problem that develops shock structures in the nonlinear Boltzmann equation.
We consider the spatial domain $x \in [0,1]$ with Neumann boundary conditions, using $N_x = 50$ points. Let the final time $t=0.2$, we test for the problem with $\varepsilon \in \{1, 10^{-4}, 10^{-8}\}$. The initial condition is discontinuous in space and defined by
\begin{equation*}
f^{\text{in}}(x, v) = 
\begin{cases}
    M_{1,\, \mathbf{0},\, 1}(v), & x \in [0, 0.5] \\
    M_{1/8,\, \mathbf{0},\, 1/4}(v), & x \in (0.5, 1]
\end{cases},
\end{equation*}
where $M_{\rho, u, T}(v)$ is the Maxwellian distribution and $\mathbf{0}=(0,0)$.

\begin{figure}
    \centering
    \begin{subfigure}[b]{0.8\textwidth}
    \caption{$\varepsilon = 10^{-4}$}
        \includegraphics[width=\textwidth]{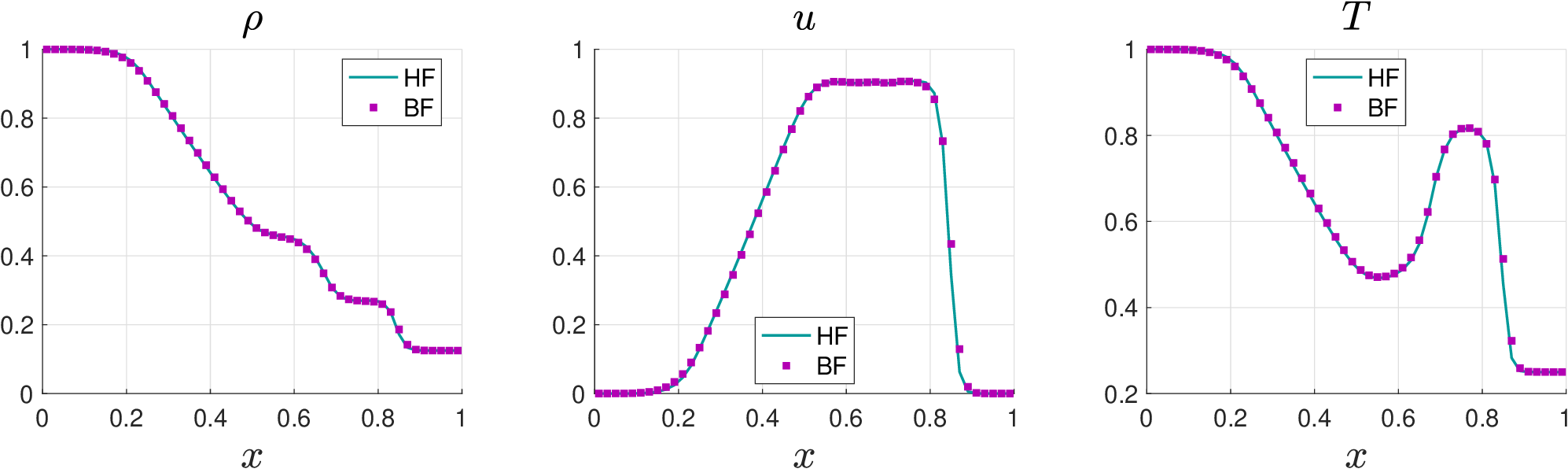}
    \end{subfigure}
    \hfill
    \hspace{1em}
    \begin{subfigure}[b]{0.8\textwidth}
    \caption{$\varepsilon = 10^{-8}$}
        \includegraphics[width=\textwidth]{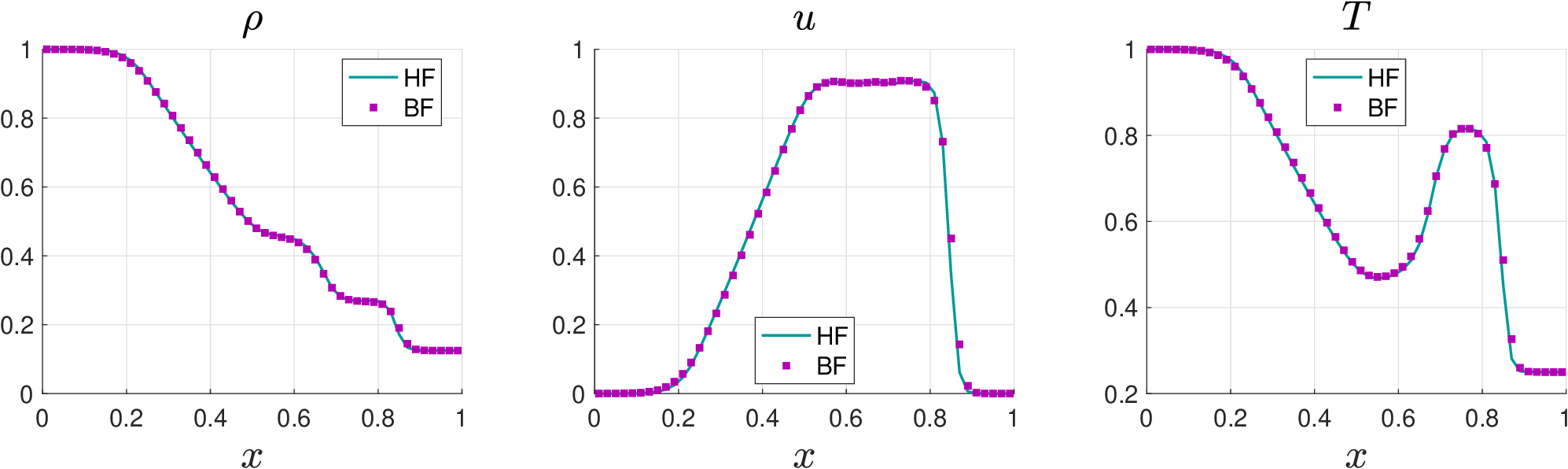}
    \end{subfigure}
    \caption{Test II: Macroscopic quantity profiles at $t = 0.2$ comparing bi-fidelity (dotted lines) and high-fidelity (solid lines) solutions for varying 
    $\varepsilon=10^{-4}$ (top row) and  $\varepsilon=10^{-8}$ (bottom row). Left: density $x\mapsto\rho(t=0.2, x)$. Middle: bulk velocity  $x\mapsto u(t=0.2, x)$. Right: temperature   $x\mapsto T(t=0.2, x)$.}
    \label{fig:Boltz_shk}
\end{figure}

Figure~\ref{fig:Boltz_shk} shows the macroscopic quantities (density, the first component of bulk velocity and temperature) at $t = 0.2$ obtained from the BF and HF solutions for varying Knudsen numbers. The BF solutions show an excellent agreement with the HF reference solution in different regimes, demonstrating the robustness of our method in solving the shock problem for the nonlinear Boltzmann equation. 

\begin{figure}
    \centering
    \begin{subfigure}[b]{0.38\textwidth}
        \includegraphics[width=\textwidth]{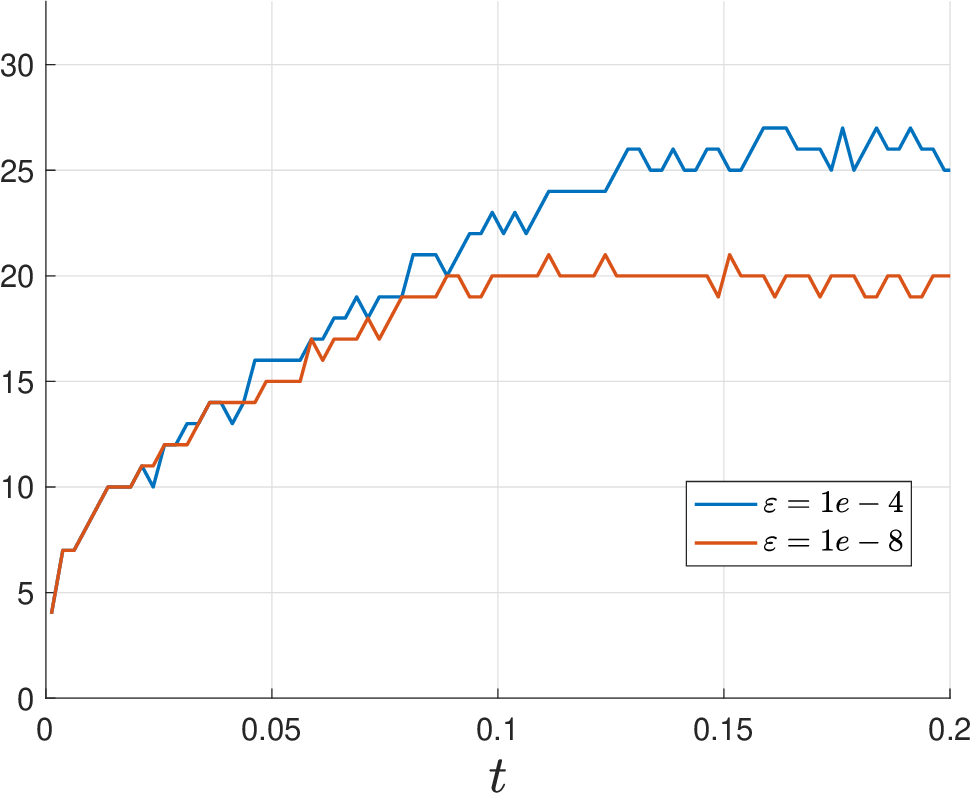}
    \end{subfigure}
    \hfill
    \begin{subfigure}[b]{0.38\textwidth}
        \includegraphics[width=\textwidth]{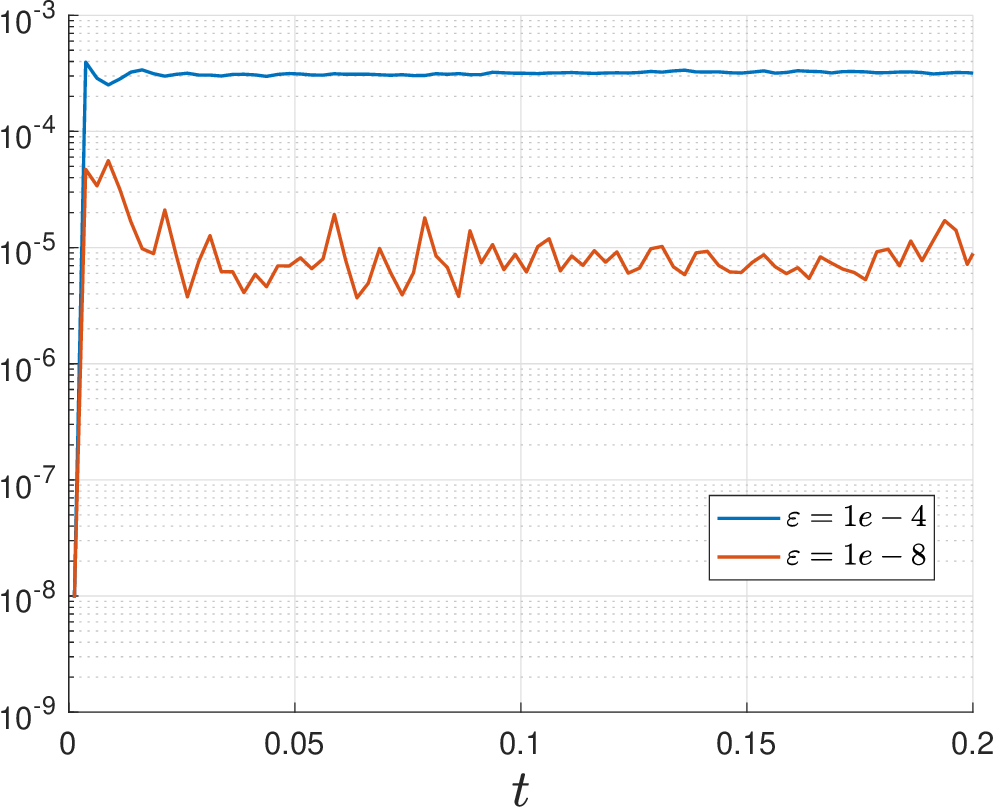}
    \end{subfigure}
    \caption{Test II: Left: Time evolution of $|\gamma|$ for different $\varepsilon$. Right: Time evolution of $\|f_B - M(f_B)\|_{\ell^1}$ for different $\varepsilon$.}
    \label{fig:Boltz_shk_history}
\end{figure}

In Figure~\ref{fig:Boltz_shk_history}, the time evolution of the size of selected velocity set $|\gamma|$ as well as the distance between BF solution and its equilibrium measure by $\|f_B - M(f_B)\|_{\ell^1}$ are plotted. As time reaches $0.1$, one can see that the smaller $\e$ is, the less points are selected in the velocity space.
Besides, we observe a gradual increase in the number of selected points as the solution converges to the equilibrium, then maintain at about $20$ or $25$.  
This is consistent with the phenomenon shown in the right-hand-side of Figure~\ref{fig:Boltz_shk_history}, that the solution is closer to the equilibrium at the initial period of time, then as time propagates the shock structure of the model begins to form, leading the BF solution stabilizes towards the macroscopic solution.

\begin{figure}
    \centering
\includegraphics[width=0.4\textwidth]{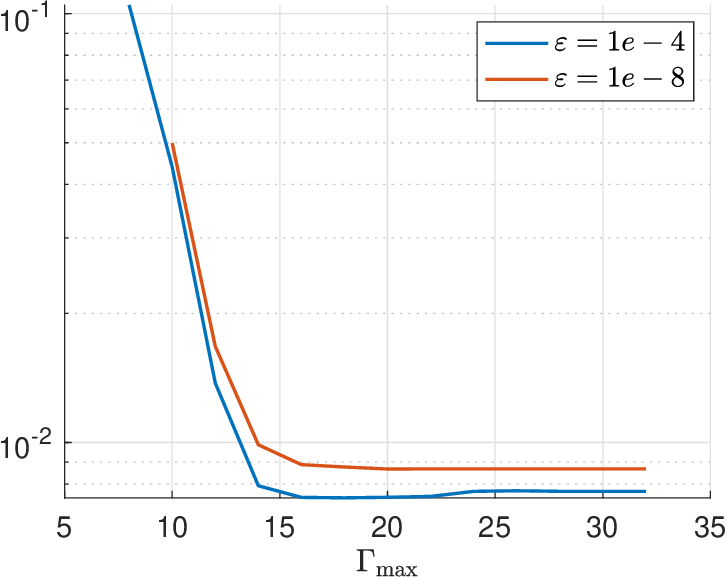}
    \caption{Test II: Relative $\mathscr{E}(t=0.2)$ error against $\Gamma_{\max}$. }
    \label{fig:Boltz_shk_error}
\end{figure}

Figure~\ref{fig:Boltz_shk_error} investigates the approximation error~\eqref{equation:error_vs_maxSelect} as the maximum allowed size $\Gamma_{\text{max}}$  increases. This plot again shows that the error between BF and HF decreases and eventually saturates as more points are included. 

\begin{figure}
    \centering
\begin{subfigure}[b]{0.3\textwidth}
    \includegraphics[width=\textwidth]{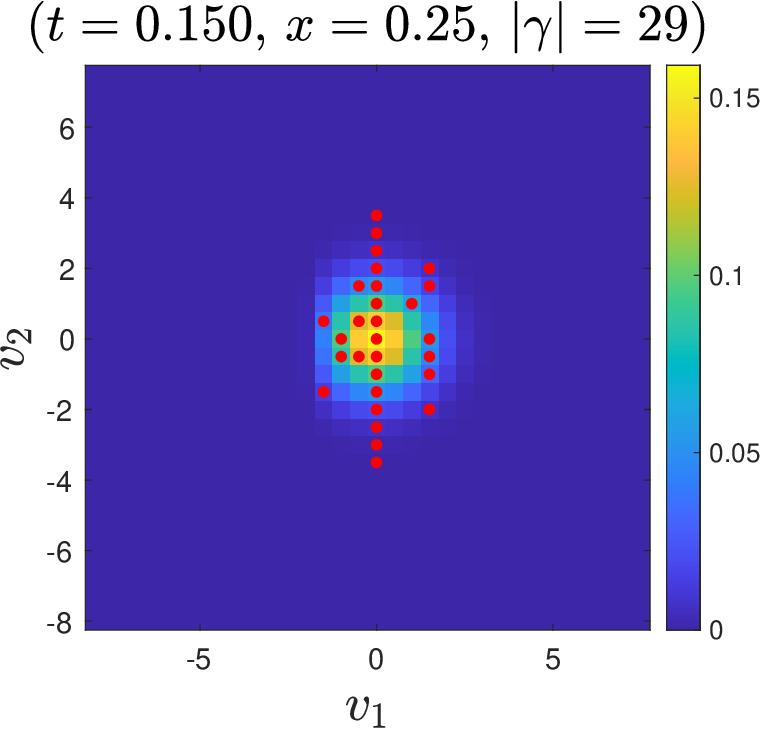}
    \end{subfigure}
    \hfill
    \hspace{1em}
    \begin{subfigure}[b]{0.3\textwidth}  
    \includegraphics[width=\textwidth]{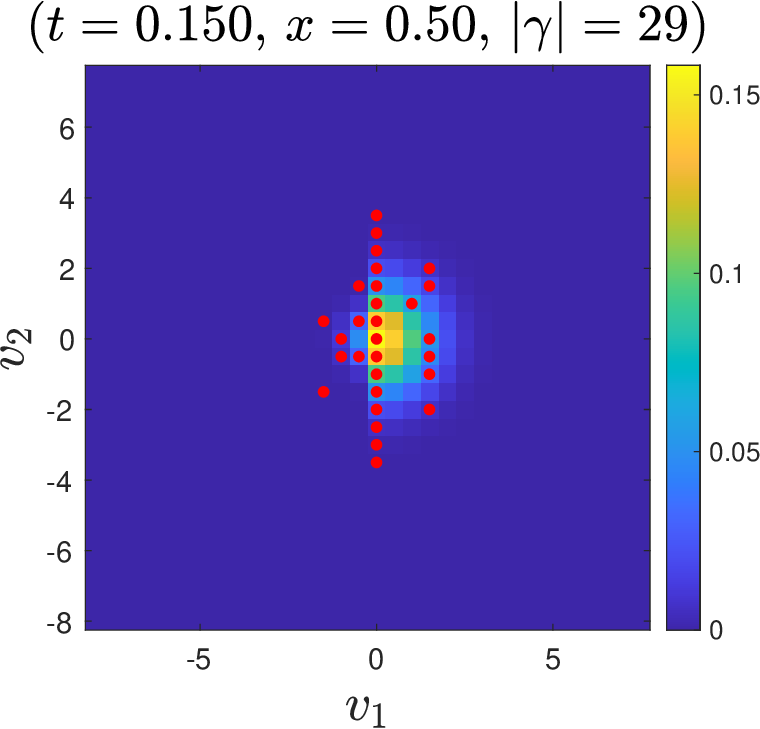}
    \end{subfigure}
    \hfill
    \hspace{1em}
    \begin{subfigure}[b]{0.3\textwidth}
    \includegraphics[width=\textwidth]{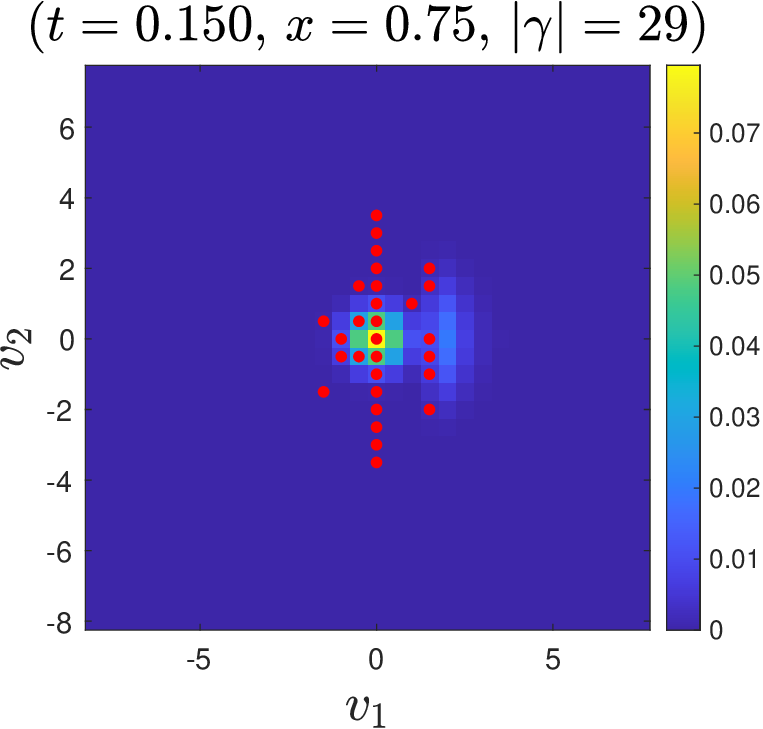}
    \end{subfigure}
    \caption{Test II: Snapshots of the layout of the selected velocity points in $\gamma$ and the shape of the distribution functions. }
    \label{fig:Gamma_dist}
\end{figure}

In Figure~\ref{fig:Gamma_dist}, we show the layout of the velocity points selected in the set $\gamma$ at three snapshots with $x\in \left\{0.25, 0.5, 0.75\right\}$ and $t=0.15$. We use the setting of Test II for $\varepsilon=10^{-8}$. The 2D color plots  are shown to represent the velocity distributions of the bi-fidelity solutions. It is noticeable from this figure that the picked velocity points can somehow capture the solution behaviour in the velocity space--more points are collected in the region where the distribution owns a more complex shape.

\subsubsection{Test III: Blast wave problem for Model 2} 

In this test, we consider a more challenging blast wave problem for the nonlinear Boltzmann equation. We set the spatial domain to be $x \in [0,1]$ and assume specular boundary conditions, using $N_x = 100$. Let the final time $t=0.1$ for $\varepsilon \in \{10^{-4}, 10^{-8}\}$. The initial condition has discontinuities in space and is given by
{\small
\begin{equation*}
f^{\text{in}}(x, v) = 
\begin{cases}
    M_{1,\, (1,0),\, 2}(v), & x \in [-0.5, -0.3] \\
    M_{1,\, (0,0),\, 1/4}(v), & x \in (-0.3, 0.3] \\
    M_{1,\, (-1,0),\, 2}(v), & x \in (0.5, 0.5] 
\end{cases},
\end{equation*}
}where $M_{\rho, u, T}(v)$ is the Maxwellian. 
 \begin{figure}
    \centering
    \begin{subfigure}[b]{0.8\textwidth}
    \caption{$\varepsilon = 10^{-4}$}
        \includegraphics[width=\textwidth]{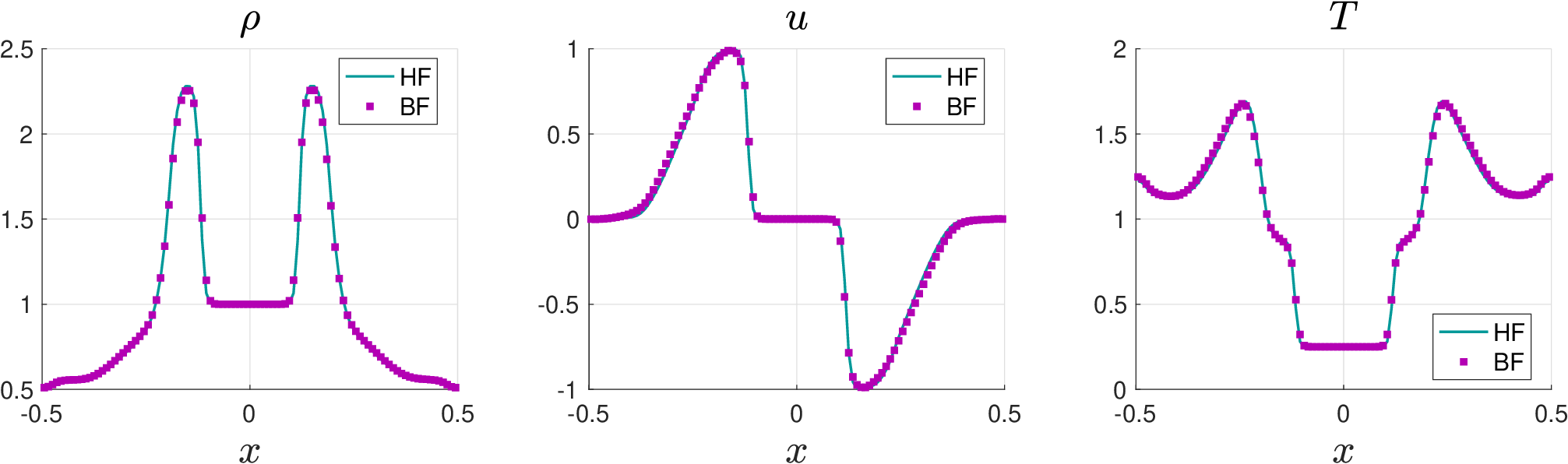}
    \end{subfigure}
    \hfill
    \hspace{1em}
    \begin{subfigure}[b]{0.8\textwidth}
    \caption{$\varepsilon = 10^{-8}$}
        \includegraphics[width=\textwidth]{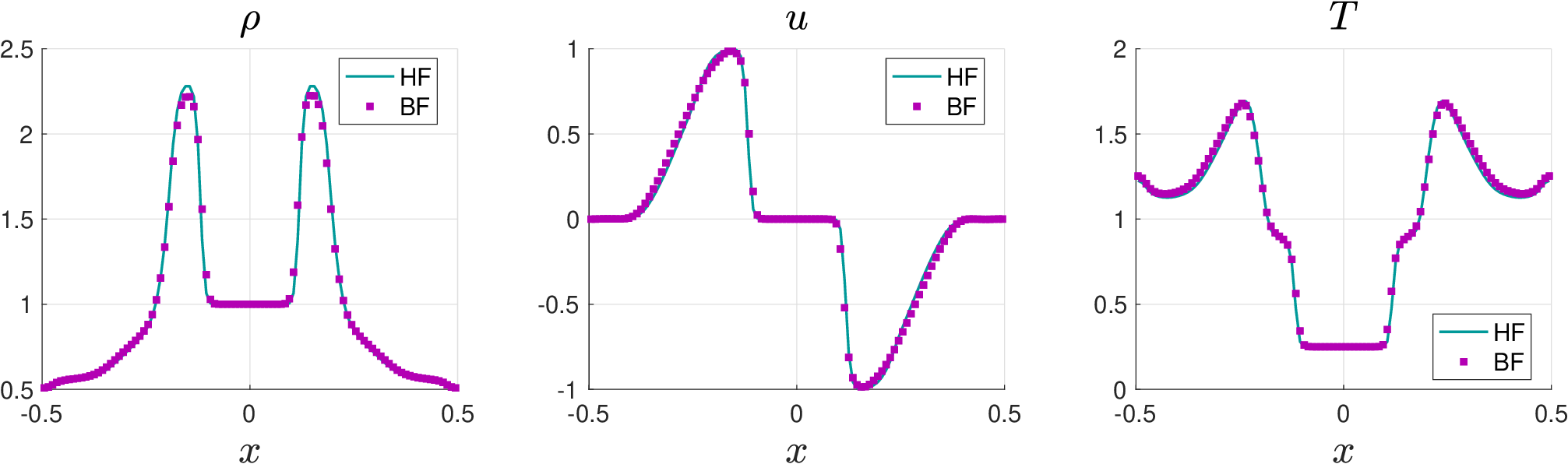}
    \end{subfigure}
    \caption{Test III: Macroscopic quantity profiles at $t = 0.1$ comparing bi-fidelity (dotted lines) and high-fidelity (solid lines) solutions for varying Knudsen numbers $\varepsilon=10^{-4}$ (top row) and  $\varepsilon=10^{-8}$ (bottom row). Left: density $x\mapsto\rho(t=0.2, x)$. Middle: bulk velocity  $x\mapsto u(t=0.2, x)$. Right: temperature   $x\mapsto T(t=0.2, x)$. }
    \label{fig:blast_wave_compare}
\end{figure}

\begin{figure}
    \centering
    \begin{subfigure}[b]{0.38\textwidth}
        \includegraphics[width=\textwidth]{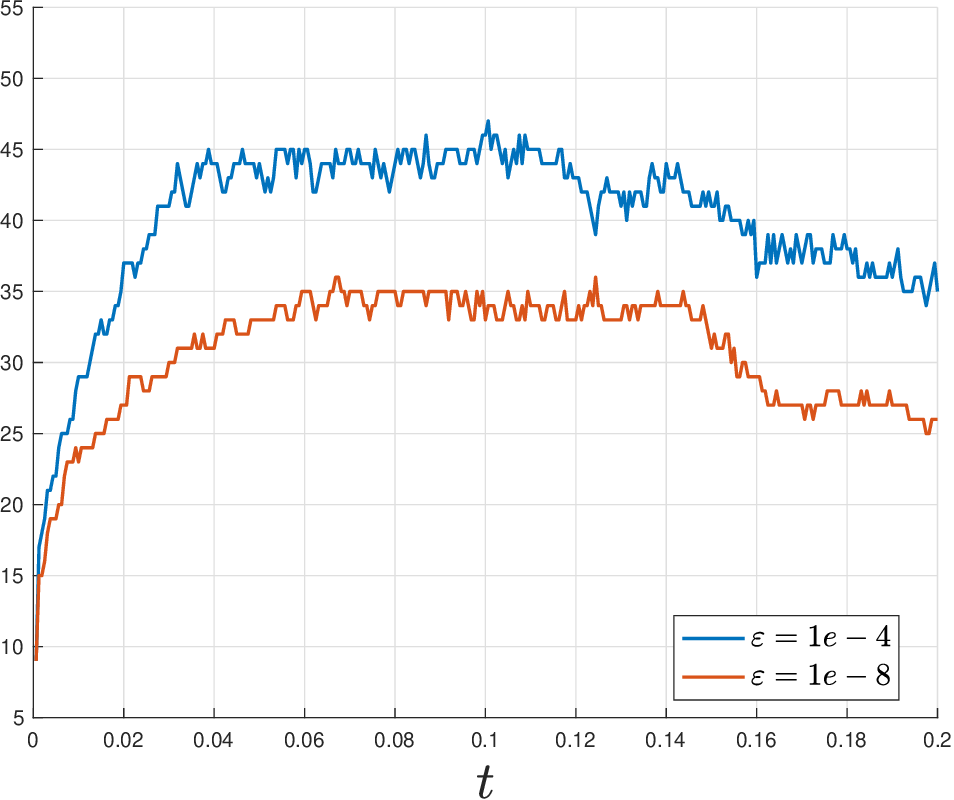}
    \end{subfigure}
    \hfill
    \begin{subfigure}[b]{0.38\textwidth}
        \includegraphics[width=\textwidth]{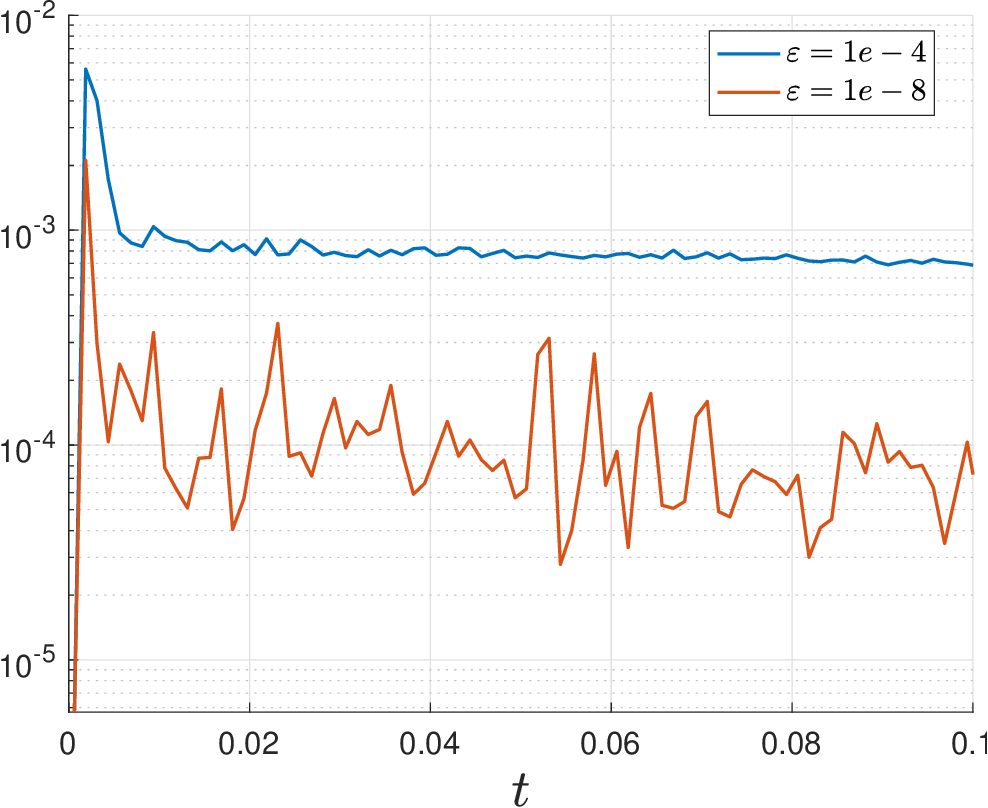}
    \end{subfigure}
\caption{Test III: Left: Time evolution of set size $|\gamma|$ for different $\varepsilon$. Right: Time evolution of $\|f_B - M(f_B)\|_{\ell^1}$ for different $\varepsilon$.}
\label{fig:Boltz_blast_history}
\end{figure}

The bi-fidelity method accurately captures the solution behaviour including the sharp shock fronts, as shown in Figure~\ref{fig:blast_wave_compare}, where the BF and HF solutions match well at all spatial points. 
Similar adaptive point selection is evident in the time evolution of $|\gamma|$ and the difference $\|f_B - M(f_B)\|_{\ell^1}$, plotted in Figures~\ref{fig:Boltz_blast_history}. The observations are similar as previous tests, which illustrates again the efficiency of our algorithm, since only $|\gamma|=45$ points are selected whereas initially the velocity grid contained $32^2=1024$ points. To the best of our efforts and based on our computing resources, we found that the BF algorithm costs only about $1/3$ compared to that of the expensive HF solver. 

\begin{figure}
    \centering
    \includegraphics[width=0.4\textwidth]{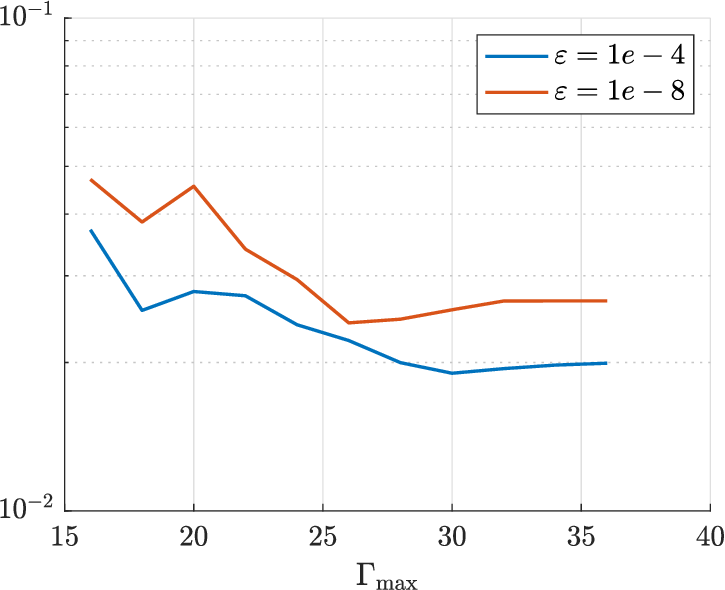}
    \caption{Test III: Relative $\mathscr{E}(t = 0.05)$ error against $\Gamma_{\text{max}}$.}
\label{fig:Boltz_blast_error}
\end{figure}

Finally, in Figure~\ref{fig:Boltz_blast_error} we provide  the approximation error~\eqref{equation:error_vs_maxSelect} as the maximum allowed size $\Gamma_{\text{max}}$ increases. Due to the complexity of this blast wave problem, the error decay is not as fast as shown in earlier tests yet still can show the effectiveness of our point selection procedure. 

\subsection{Empirical Error Estimation}

In this section, we investigate the empirical error estimation in section~\ref{section:empirical_error_estimation} as a metric to further demonstrate the accuracy of bi-fidelity method. We study the smooth initial problem for the nonlinear Boltzmann equation, i.e., Test I (b). Let the final time $t=0.1$, and $\varepsilon \in \{1, 10^{-4}, 10^{-8}\}$. The time evolutions of the following four quantities are computed: (i) the similarity ratio $R_s$ defined in~\eqref{equation:rs}; (ii) the projective distance $\max_v R_e(v)$ defined in~\eqref{equation:re}; (iii) the empirical error bound $\mathcal{E}$ defined in~\eqref{equation:empirical_error}; (iv) The `true error' $\|f_B - f_H\|_{\ell^2_x} / \|f_H\|_{\ell^2_x}$ between the bi-fidelity and high-fidelity updates. We will plot these quantities as time evolves. 

\begin{figure}
    \centering
    \begin{subfigure}[b]{0.3\textwidth}        \includegraphics[width=\textwidth]{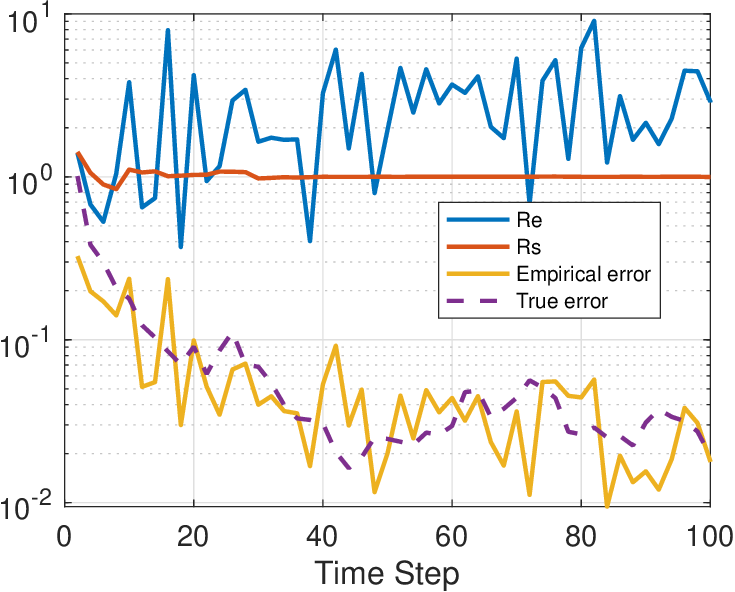}
        \caption{$\varepsilon = 1$}
    \end{subfigure}
    \hfill
    \begin{subfigure}[b]{0.3\textwidth}
\includegraphics[width=\textwidth]{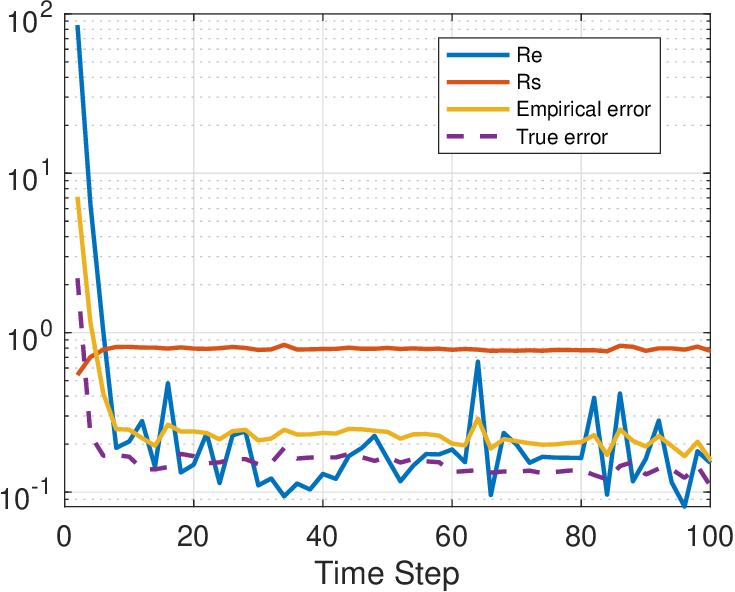}
        \caption{$\varepsilon = 10^{-4}$}
    \end{subfigure}
    \hfill
    \begin{subfigure}[b]{0.3\textwidth}
\includegraphics[width=\textwidth]{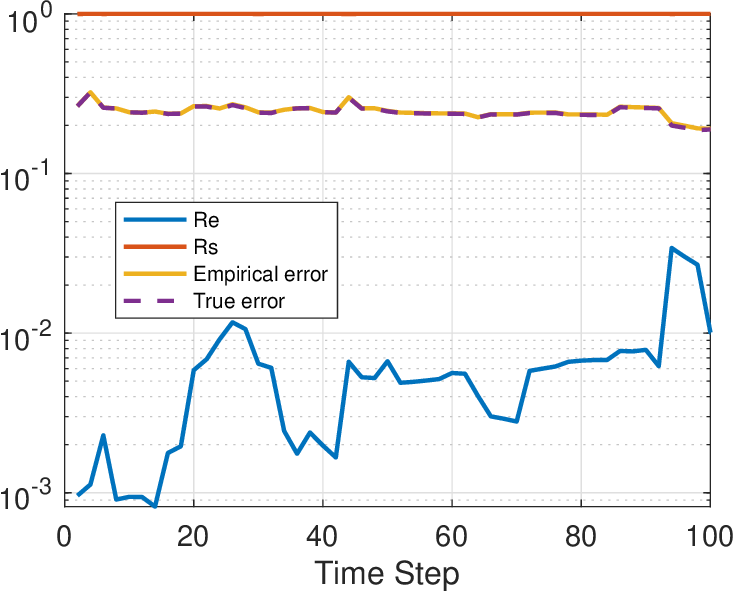}
        \caption{$\varepsilon = 10^{-8}$}
    \end{subfigure}
    \caption{Empirical error bound test: time evolutions of $R_s, \max_v R_e(v)$, empirical error $\mathcal{E}$ and `true error' $\|f_B - f_H\|_{\ell^2_x} / \|f_H\|_{\ell^2_x}$. From left to right: $\varepsilon = 1, 10^{-4}, 10^{-8}$.}    \label{fig:empirical_error_smooth}
\end{figure}

In Figure~\ref{fig:empirical_error_smooth}, we observe that the similarity ratio $R_s$ remain close to $1$ for different $\varepsilon$ cases, which guarantees that the error bound~\eqref{equation:empirical_error_bound} serves as an effective and robust surrogate of the true error. 
Additionally, $\max_v R_e(v)$ is bounded by $10$, indicating that the in-plane error does not dominate over the projective distance, and the BF solution is expected to deliver satisfactory results. 
Figure~\ref{fig:empirical_error_smooth} also shows us that  the empirical error bound $\mathcal{E}$ indeed evolves in close proximity to the true error as time step increases, further demonstrating that the empirical error bound is a reliable surrogate for the true error. One observes that as $\varepsilon$ decreases, the magnitude of $\max_v R_e(v)$ decreases correspondingly, resulting in a more accurate empirical error bound $\mathcal{E}$. We remark that this empirical error bound $\mathcal{E}$ requires no additional computation in the bi-fidelity algorithm.

%% file: 6_Conclusion.tex
\section{Conclusion}
In this work, we present a bi-fidelity strategy for velocity discretization of Boltzmann-type equations, which can obtain the numerical solutions efficiently and accurately.
Through the pivoted Cholesky decomposition, our chosen low-fidelity BGK model is capable of identifying the significant velocity points. Then we use the asymptotic-preserving scheme and numerically compute the updates on these few points in the velocity space, yielding accurate reconstructions at a much lower computational cost. The asymptotic-preserving analysis and empirical error bounds are shown, which provide a computable error estimates and ensure the numerical scheme being asymptotic-preserving. A series of numerical experiments on linear semiconductor and nonlinear Boltzmann problems with smooth or discontinuous initial conditions and under various regimes have been carefully studied, which demonstrates the effectiveness and robustness of our proposed method.

In the future work, we will explore different ways of point selection and study higher dimensional problems.

%% file: 7_Appendix.tex
\appendix

\section{A review of previous bi-fidelity method for UQ problems}
\label{appendix:UQ}

In the framework of stochastic collocation, the non-intrusive multi-fidelity method was studied in~\cite{NarayanGittelsonXiu2014} for general PDEs containing uncertainties. 
Below is a review of the bi-fidelity algorithm in the UQ problem setting.
Consider a time-dependent differential equation that contains uncertain parameters: 
\begin{equation*}
\begin{cases}
\mathbf{u}_t(x, t, z) = \mathcal{L}(\mathbf{u}), & \text{in } D \times (0, T] \times I_z, \\
\mathbf{u} = \mathbf{u}_0, & \text{in } D \times \{ t = 0 \} \times I_z,
\end{cases}
\end{equation*}
where $\mathcal{L}$ is a general operator that can be linear or nonlinear, 
$\mathbf{u}_0$ is the initial condition, and $z$ is some finite-dimensional random variable obtained from parametrization of the stochastic process \cite{KL}. 
The quantity of interest is the solution $\mathbf{u}(z)$ at final time. The bi-fidelity solution at final time is approximated by
{\small
\begin{equation*}
\mathbf{u}_B(z) = \sum_{k=1}^{K} c_k^L(z) \mathbf{u}_H(z_k),
\end{equation*}
}where the dependence of $t$ and $x$ in the solution is omitted. Here $\mathbf{u}_H$ is the high-fidelity (goal) model, $\mathbf{u}_B$ is the bi-fidelity approximation and $c_k^L(z)$ is the projection coefficient determined by the low-fidelity solution $\mathbf{u}_L$. 
In the offline stage, one employs the cheaper low-fidelity model--which is chosen to our own favor and problem dependent--to select the most important parameter points $\{\gamma\}$ by the greedy procedure~\cite{DeVore-Petrova-Wojtaszczyk-2013}. In the online stage, for any point $z \in I_z$, we project the low-fidelity solution $\mathbf{u}_L(z)$ onto the low-fidelity approximation space $U_L(\gamma):=\{ \mathbf{u}_L(z_1), \cdots, \mathbf{u}_L(z_K)\}$, namely
{\small
\begin{equation*}
\mathbf{u}_L(z) \approx P_{U_L(\gamma)}[\mathbf{u}_L(z)] = \sum_{k=1}^{K} c_k^L(z) \mathbf{u}_L(z_k),
\end{equation*}
}where $P_{U_L(\gamma)}$ is the projection operator onto $U_L(\gamma)$ and $\{c_k^L\}$ are the projection coefficients  computed by the Galerkin approach: $\mathbf{G}^L \mathbf{c} = \mathbf{f}$, where $\mathbf{f} = (f_k)_{1 \leq k \leq K}$ with $\displaystyle f_k = \langle \mathbf{u}_L(z), \mathbf{u}_L(z_k) \rangle$ and the Gramian matrix $\mathbf{G}^L$. We refer details to earlier work \cite{LPZ2022, BLPZ2022, Liu-Zhu-2020} on multiscale kinetic problems with uncertainties. 
We emphasize that in these previous work on UQ problems,  the bi-fidelity method is employed in the {\it random space} discretization and one requires to select the important points of random variable {\it only once}, at the final time. 

\section{The DVM method for discretization of $\mathcal{Q}_{\text{NB}}(f, f)$} 
\label{appendix:DVM}
We review the DVM method~\cite{Panferov-Heintz-2002} for the nonlinear Boltzmann collision operator
$\mathcal{Q}_{\text{NB}}(f, f)$
for $d_v \geq 2$.
The DVM is based on the Carleman representation~\cite{Carleman1932}. By defining new variables $u = v^\prime  - v$ and $w = v^\prime_\star - v$, we obtain the following form
$$
\mathcal{Q}_{\text{NB}}(f, f)(v) = \int_{\mathbb{R}^{d_v}} \int_{E_u} B^c(u, w) \left( f(v + u) f(v + w) - f(v) f(v + u + w) \right)\mathrm{d}w \mathrm{d}u 
$$
where
$
B^c(u, w) = 2|u|^{-2} B\Big(\sqrt{u^2 + w^2}, \frac{|u|}{\sqrt{u^2 + w^2}}\Big)
$
and $E_u$ is the plane through the origin perpendicular to the vector $u$.

We denote 
$
F(v, u, w) = B^c(u, w) \left( f(v + u) f(v + w) - f(v) f(v + u + w) \right)
$
and 
\begin{equation}
    \label{eq:G}
    G(v, u) = \int_{E_u} F(v, u, w) \mathrm{d} w.    
\end{equation}
The collision operator $\displaystyle \mathcal{Q}_{\text{NB}}(f, f)(v) = \int_{\mathbb{R}^{d_v}} G(v, u) \mathrm{d} u $ is then approximated by 
$$
\mathcal{Q}_{\text{NB}}(f, f)(v_k) \approx (\Delta v)^{d_v} \sum_{i \in \mathbb{Z}^{d_v}} G(v_k, u_i).
$$
For fixed $v_k$ and $u_i$, the integral~\eqref{eq:G} is approximated by 
$$
G(v_k, u_i) \approx (\Delta v)^{d_v - 1} \mathrm{det}(L_i) \sum_{j \in L_i} F(v_k, u_i, w_j),
$$
where $L_i = \{j\in\mathbb{Z}^{d_v} \mid i \cdot j = 0\}$ is the lattice that is the intersection of the velocity grid with the plane $E_i$. We have an explicit expression
$\mathrm{det}(L_i) = \frac{|i|}{\mathfrak{g}(i)}$ where $\mathfrak{g}(i)$ is the greatest common divisor of the components of $i$. 
We obtain the formal approximation: for $k\in\mathcal{N}$,
{\small
\begin{equation}
\label{equation:DVM}
\begin{aligned}
        \mathcal{Q}_{\text{NB}}(f, f)(v_k) \approx &  (\Delta v)^{2d_v - 1} \sum_{i\in \mathbb{Z}^{d_v}} \mathrm{det} (L_i) \sum_{j\in L_i} B^c_{i, j} \\
        &\times\Big(f(v_k + u_i)f(v_k + w_j) - f(v_k)f(v_k + u_i + w_j)\Big), 
\end{aligned}
\end{equation}
}with $B^c_{i, j} = B^c(u_i, w_j)$.
In the bi-fidelity Algorithm \ref{alg:bifid_nonlinearBoltzmann}, one only needs to evaluate the complicated integration~\eqref{equation:DVM} for $v_k\in\gamma^{n+1}$, whereas the summation on the right hand side is for $i\in\mathcal{N}$ which ensures the conservation properties of the collision operator.

\section{The AP property of the bi-fidelity algorithm for Model 1}
\label{appendix:AP_Model1}
We prove the weak AP property of the Algorithm~\ref{alg:bifid_semi} for Model 1.
\begin{theorem}
\label{theorem:AP_Model1}
Suppose $\delta < \varepsilon$
and let $(f_B^n)_n$ be the numerical solution obtained from Algorithm 
\ref{alg:bifid_semi}. Then, $(f^n_B)_n$ satisfies the following property:
$$
\left|f^n_B - \rho(f^n_B)M \right| = \mathcal{O}(\varepsilon) \Rightarrow
\left|f^{n+1}_B - \rho(f^{n+1}_B)M \right| = \mathcal{O}(\varepsilon). 
$$
In particular, as $\varepsilon \rightarrow 0$, 
our numerical scheme automatically becomes a consistent discretization of the macroscopic limit equation, thus satisfies the AP property. 

\begin{proof}

From the even-odd decomposition \eqref{eq:evenodd} we know $\rho(r_B) = \rho(f_B)$.
By the assumption $f^n_B - \rho(f^n_B)M = \mathcal{O}(\varepsilon)$, we have $j^n_B = \mathcal{O}(1)$ and $r^n_B = \rho(r^n_B)M + \mathcal{O}(\varepsilon)$.
Using the first equation in \eqref{scheme:highfid_semi} and
following a similar argument as in Theorem~\ref{theorem:AP}, we can show that  
\begin{equation}
\label{eq:r_AP}
    r^{n+1}_B = \rho(r^{n+1}_B)M + \mathcal{O}(\varepsilon + \delta) = \rho(f^{n+1}_B)M + \mathcal{O}(\varepsilon + \delta).
\end{equation}
Since $\phi = 1$, using the relationship~\eqref{eq:r_AP} in the scheme~\eqref{scheme:highfid_semi} for $j^{n+1}_B$, we deduce that $j^{n+1}_B = \mathcal{O}(1)$.
Finally, since $f^{n+1}_B = r^{n+1}_B + \varepsilon j^{n+1}_B$, we have
\begin{equation}
    f^{n+1}_B - \rho(f^{n+1}_B)M = \mathcal{O}(\varepsilon + \delta).
\end{equation}
By choosing $\delta = \mathcal{O}(\varepsilon)$, we have $f^{n+1}_B - \rho(f^{n+1}_B)M = \mathcal{O}(\varepsilon)$.

\end{proof}
\end{theorem}